\documentclass[12pt,notitlepage]{amsart}
\usepackage{latexsym,amsfonts,amssymb,amsmath,amsthm} 
\newcommand{\mz}{\ensuremath{\mathbb Z}}
\newcommand{\mr}{\ensuremath{\mathbb R}}
\newcommand{\mh}{\ensuremath{\mathbb H}}

\newcommand{\mc}{\ensuremath{\mathbb C}}
\newcommand{\mf}{\ensuremath{\mathbb F}}

\newcommand{\shortmod}{\ensuremath{\negthickspace \negthickspace \negthickspace \pmod}}

\newcommand{\half}{\ensuremath{ \frac{1}{2}}}

\newcommand{\intR}{\int_{-\infty}^{\infty}}
\newcommand{\sign}{\mathfrak{a}}
\newcommand{\chibar}{\overline{\chi}}
\newcommand{\sumstar}{\sideset{}{^*}\sum}

\theoremstyle{plain}		
	\newtheorem{mytheo}{Theorem}[section]
	
	\newtheorem{myprop}[mytheo]{Proposition}
	
     \newtheorem{mylemma}[mytheo]{Lemma}
	\newtheorem{mydefi}[mytheo]{Definition}
	\newtheorem{myconj}[mytheo]{Conjecture}
	
\theoremstyle{remark}

\numberwithin{equation}{section}
\begin{document}
\today
\title{The fourth moment of Dirichlet L-functions}
\author{Matthew P. Young} 
\address{American Institute of Mathematics, 360 Portage Ave.,
Palo Alto, CA 94306-2244}
\curraddr{Department of Mathematics \\
 	  Texas A\&M University \\
 	  College Station \\
	  TX 77843-3368 \\
		U.S.A.}
\email{myoung@math.tamu.edu}
\thanks{This research was supported by an NSF Mathematical Sciences Post-Doctoral Fellowship and by the American Institute of Mathematics.}
\begin{abstract}
We compute the fourth moment of Dirichlet L-functions at the central point for prime moduli, with a power savings in the error term.
\end{abstract}
\maketitle
\section{Introduction}
Estimating moments of families of L-functions is a central problem in number theory due in large part to extensive applications.  Yet, these moments are seen to be natural objects to study in their own right as they illuminate structure of the family and display beautiful symmetries. 

The Riemann zeta function has by far garnered the most attention from researchers.  Ingham \cite{Ingham} proved the asymptotic formula
\begin{equation*}
\frac{1}{T} \int_0^{T} |\zeta({\textstyle \half + it})|^4 dt = a_4 (\log{T})^4 + O((\log{T})^3)
\end{equation*}
where $a_4 = (2 \pi^{2})^{-1}$.  Heath-Brown \cite{H-B2} improved this result by obtaining
\begin{equation}
\label{eq:zetapowersavings}
\frac{1}{T} \int_0^{T} |\zeta({\textstyle \half + it})|^4 dt = \sum_{i=0}^{4} a_i (\log{T})^i + O(T^{-\frac{1}{8} +\varepsilon})
\end{equation}
for certain explicitly computable constants $a_i$ (see (5.1.4) of \cite{CFKRS}).
Obtaining a power savings in the error term was a significant challenge because it requires a difficult analysis of off-diagonal terms which contribute to the lower-order terms in the asymptotic formula.  

The fourth moment problem is related to the problem of estimating
\begin{equation}
\label{eq:dust}
\sum_{n \leq x} d(n) d(n+f)
\end{equation}
uniformly for $f$ as large as possible with respect to $x$.  Extensive discussion of this binary divisor problem can be found in \cite{Motohashi1}.
The sum \eqref{eq:dust} can be transformed into a problem involving Kloosterman sums.  The strength of Heath-Brown's result \eqref{eq:zetapowersavings} relies on the Weil bound. 
Using the spectral theory of automorphic forms (the Kuznetsov formula), Iwaniec (\cite{Iwaniec}, Theorem 3) showed
\begin{equation}
\label{eq:zetashortinterval}
\int_{T}^{T + T^{2/3}} |\zeta({\textstyle \half + it})|^4 dt \ll T^{2/3 + \varepsilon}.
\end{equation}
Notice that Heath-Brown's result gives \eqref{eq:zetashortinterval} with $T^{2/3}$ replaced by $T^{7/8}$.  With some extra work, these mean-value estimates can be turned into subconvexity estimates for the zeta function (although Weyl's method already gives $\zeta({\textstyle{\half + it}}) \ll t^{1/6 + \varepsilon}$ with far less effort).

Motohashi has proved a beautiful exact formula for a smoothed version of the fourth moment of zeta in terms of the cubes of the central values of all automorphic (degree $2$, level $1$) L-functions (\cite{M}, Theorem 4.2).  The result is difficult to describe precisely because the dependence on the smoothing function involves elaborate integral transforms.  Motohashi considers
\begin{equation}
\label{eq:motohashishift}
M_g(\alpha,\beta,\gamma,\delta) = \int_{-\infty}^{\infty} \zeta({\textstyle \half + \alpha + it}) \zeta({\textstyle \half + \beta + it})\zeta({\textstyle \half + \gamma - it})\zeta({\textstyle \half + \delta - it}) g(t/T) dt,
\end{equation}
where initially the parameters $\alpha,\beta,\gamma,\delta$ have large real parts, and uses the Kuznetsov formula to develop a meromorphic continuation of $M_g(\alpha,\beta,\gamma,\delta)$ to include the origin.  These parameters (or shifts) are helpful for developing the main term in the asymptotic formula for the fourth moment because certain residue computations become simplified (e.g., it is much easier to compute the residues of five simple poles than it is to compute a residue at a fifth order pole).  The smoothing function $g$ can be replaced with the sharp truncation $t \leq T$ using standard techniques (however, unsurprisingly, the error term becomes larger).

Recently, five authors \cite{CFKRS} have developed a recipe for conjecturing moments of families of L-functions, including all lower order terms.  The presence of these shifts plays an important role in their conjectures.  Besides allowing for simpler residue computations, the presence of these parameters makes symmetries of the family visible.  For example, $M_g$ is invariant under switching $\alpha$ and $\beta$ or $\gamma$ and $\delta$.  Application of the functional equation for zeta shows relations under $\alpha \rightarrow -\gamma$, $\gamma \rightarrow -\alpha$, and other similar relations.  The shifts are useful too since for example one can differentiate with respect to them to study moments of the derivatives of a family.

The family of all primitive Dirichlet L-functions of modulus $q$ is analogous in some ways to the Riemann zeta function in $t$-aspect.  However, there are significant differences which cause the family of Dirichlet L-functions to be more difficult to study.  Ramachandra \cite{R} conjectured that for prime $p$,
\begin{equation*}
p^{-1} \sum_{\chi \shortmod{p}} |L({\textstyle \half },\chi)|^4 \sim (2\pi^2)^{-1} (\log{p})^4.
\end{equation*}
Heath-Brown \cite{H-B1} proved
\begin{equation}
\label{eq:heathbrown}
\frac{1}{\phi^*(q)} \sumstar_{\chi \shortmod{q} }|L({\textstyle \half },\chi)|^4 = \frac{1}{2\pi^2}
\prod_{p | q} \frac{(1-p^{-1})^3}{(1+p^{-1})}
 (\log{q})^4 + O\Big(2^{\omega(q)} \frac{q}{\phi^*(q)} (\log{q})^3\Big),
\end{equation}
where $\phi^*(q)$ is the number of primitive characters modulo $q$, $\omega(q)$ is the number of distinct prime factors of $q$, and the sum is over all primitive characters modulo $q$.  For almost all $q$ this is an asymptotic formula, however if $\omega(q) \geq (\log \log{q})/\log{2}$ the error term is not smaller than the main term.  Recently, Soundararajan \cite{S} improved the error term in \eqref{eq:heathbrown} so that an asymptotic formula does indeed hold for all $q$.

In this paper we obtain the asymptotic formula with a power savings for prime moduli.
\begin{mytheo} 
\label{thm:mainresult}
For prime $q \neq 2$, we have
\begin{equation}
\label{eq:mainresult}
\frac{1}{\phi^*(q)} \sumstar_{\chi \shortmod{q} }|L({\textstyle\frac 1{2}},\chi)|^4 = \sum_{i=0}^{4} c_i (\log{q})^i + O(q^{-\frac{5}{512} + \varepsilon}),
\end{equation}
for certain computable absolute constants $c_i$.  Here $-5/512 = -1/80 + \theta/40$, where $\theta = 7/64$ is the current best-known bound on the size of the Hecke eigenvalue $\lambda(n)$ of a Maass form, that is, a bound of the form $|\lambda(n)| \leq d(n) n^{\theta}$.
\end{mytheo}
This result is the analog of \eqref{eq:zetapowersavings} for Dirichlet L-functions of prime moduli.  The method of proof can possibly generalize to handle certain variations on \eqref{eq:mainresult} with twists of fixed Dirichlet L-functions.  Such variations could be used to prove simultaneous nonvanishing results for four Dirichlet L-functions.  For this type of variation it is imperative to have a power saving in the error term because the main term may be as small as a constant.

The study of \eqref{eq:mainresult} with $L(1/2, \chi)^2$ replaced by $L(f \times \chi, 1/2)$ for $f$ a Hecke cusp form has similarities to the analysis in this paper.  The difference amounts to replacing the divisor function with the Hecke eigenvalues of the cusp form.  In this paper we crucially use the fact that the divisor function can be written as $1\star 1$.  It is a difficult and interesting problem to obtain a power savings in the second moment of Dirichlet twists of a Hecke cusp form.

In order to obtain power savings in \eqref{eq:mainresult} it is necessary to obtain an asymptotic formula for the off-diagonal terms.  The results of Heath-Brown and Soundararjan arise from bounding the contribution of these off-diagonal terms.  To elaborate, the problem of estimating \eqref{eq:mainresult} essentially reduces to the analysis of the following divisor sum
\begin{equation}
\label{eq:divisorsumexample}
\frac{1}{\phi(q)} \sum_{d | q} \phi(d) \mu(q/d) \mathop{\sum \sum}_{\substack{m \equiv n \shortmod{d} \\ (mn,q) = 1}} \frac{d(m) d(n)}{\sqrt{mn}} V\left(\frac{mn}{q^2}\right),
\end{equation}
where $V(x)$ is a smooth function with rapid decay, satisfying $V(x) \approx 1$ for $x$ small.  
The contribution of the diagonal terms $m=n$ in the sum is easily computed and gives the main term of \eqref{eq:heathbrown}.

The off-diagonal terms $m \neq n$ are much more difficult to analyze.  One of the primary difficulties in treating \eqref{eq:divisorsumexample} is the large ranges of summation of $m$ and $n$.  Consider the sum $B_{M,N}$ consisting of those summands in \eqref{eq:divisorsumexample} in the dyadic segments $M < m \leq 2M$, $N < n \leq 2N$, where $MN \leq q^2$ (by symmetry suppose $N \geq M$).  If $M \asymp 1$ then $N$ can be as large as $q^2$ but if $M \asymp q$ then $N \asymp q$ also.  It is not surprising that different techniques are required to estimate $B_{M,N}$ for these different regions.  In fact, the relative sizes of $M$ and $N$ has a large effect on how large $B_{M,N}$ is.  Qualitatively, $B_{M,N}$ is small if $M$ and $N$ are far from each other, but when $M$ and $N$ are close to the same size, then $B_{M,N}$ contributes to the main term of \eqref{eq:mainresult}.  To see this former assertion heuristically, we use the fact that the divisor function is uniformly distributed in arithmetic progressions, that is if $(m,q)=1$ 
\begin{equation}
\label{eq:divisorsumAP}
\sum_{\substack{n \leq x \\ n \equiv m \shortmod{q}}} d(n) = \frac{1}{\phi(q)} \sum_{\substack{n \leq x \\ (n,q) = 1}} d(n) + (\text{error}).
\end{equation}
It is known unconditionally that \eqref{eq:divisorsumAP} holds for $q < x^{\frac{2}{3} - \varepsilon}$ with an error of size $O(x^{\frac13 + \varepsilon})$, following from Weil's bound for Kloosterman sums; see Corollary 1 of \cite{H-B2} for a proof.  Note that if $q$ is prime, then
\begin{equation*}
B_{M,N} \approx \frac{1}{\sqrt{MN}} \sum_{\substack{M < m \leq 2M \\(m, q) = 1 }} d(m)
\Big(
 \mathop{\sum}_{\substack{N < n \leq 2N \\ n \equiv m \shortmod{q}} }  d(n) 
 - \frac{1}{\phi(q)} \sum_{\substack{N < n \leq 2N \\ (n, q) = 1 } } d(n)
\Big),
\end{equation*} 
which is small provided the error term of \eqref{eq:divisorsumAP} is smaller than the main term.

Improving the uniformity for which \eqref{eq:divisorsumAP} is true is a challenging open problem.  It seems natural to expect that the divisor function is evenly distributed across arithmetic progressions for $q < x^{1 - \varepsilon}$.  Fouvry \cite{Fouvry} proved this is true on average over $x^{2/3 + \varepsilon} < q < x^{1-\varepsilon}$ , and Fouvry and Iwaniec \cite{FI} have produced results which cover $x^{2/3 - \varepsilon} < q < x^{2/3 + \varepsilon}$ for special values of $q$ .  These results fix the arithmetic progression $m \pmod{q}$ and average over $q$; for the application of estimating \eqref{eq:divisorsumexample}, $q$ is fixed but $m$ is allowed to vary.  
One of the main difficulties in this work is treating the range of summation in \eqref{eq:divisorsumexample} where $M \asymp q^{\half}$, $N \asymp q^{\frac32}$, which is at the edge of the range where the error term of $O(x^{\frac13 + \varepsilon})$ for \eqref{eq:divisorsumAP} is barely insufficient.

The expectation \eqref{eq:divisorsumAP} breaks down when $x$ and $q$ are of comparable size.  Instead, main terms are formed from the summands with $M \approx N \approx q$.  To see this heuristically, note that in this range the sum \eqref{eq:divisorsumexample} is mimicked by
\begin{equation*}
\sumstar_{a \shortmod{q}} \Big( \sum_{\substack{m \approx q \\m \equiv a \shortmod{q}}} \frac{d(m)}{\sqrt{m}} - \frac{1}{\phi(q)} \sum_{m \approx q} \frac{d(m)}{\sqrt{m}} \Big)^2,
\end{equation*}
which can be seen by opening the square and summing over $a$.
Clearly this expression cannot be small because the sum of $m \approx q$, $m \equiv a \pmod{q}$ essentially picks up one term which certainly cannot approximate the average behavior of the divisor function!

This dichotomy is somewhat analogous to the (smoothed) fourth moment of the zeta function, where averaging over $t$ aspect forces $m$ and $n$ to be close to each other (because there is an integral of the form $\int (n/m)^{it} g(t) dt$, which is small unless $m$ and $n$ are close).  It is remarkable that the same phenomenon occurs for Dirichlet L-functions (the main contribution coming from $m$ and $n$ close) but for completely different reasons.  If we combined our average over characters with a short $t$-aspect integration then our problem would become much simpler; however, this would destroy any nonvanishing applications.  

One apparent difficulty in finding an asymptotic formula for \eqref{eq:divisorsumexample} is that the values of $V(mn/q^2)$ when $m$ and $n$ are around $q$ are at the transition range of $V$.  However, this type of behavior has been encountered by Soundararajan (\cite{S2}, discussion following (5.16)) and Kowalski (\cite{Kowalski}, p. 155) in their studies of other families of L-functions.  The analysis of the main term in this region of summation leads to a contour integral which can be computed exactly using symmetry properties of the integrand, which in turn relies upon the functional equation for the Riemann zeta function.  In our work this occurs in the proof of Lemma \ref{lemma:QandQminus}.

In this paper we use essentially two different methods of handling the sum \eqref{eq:divisorsumexample}, depending on the sizes of $M$ and $N$.  When $N$ is significantly larger than $M$, then we treat the sum over $n$ as the divisor function in arithmetic progression.  Using only Weil's bound in a straightforward way, we could obtain the necessary asymptotic formula for $N/M > q^{1 + \varepsilon}$.  We succeeded at extending the range of summation of $n$ to smaller values by exploiting extra savings by averaging over $m$.  These arguments are presented in Section \ref{section:divisorAP}.  

The region where $M$ and $N$ are relatively close is treated using similar techinques as Motohashi in his work on the fourth moment of the Riemann zeta function \cite{M}.  The Kuznetsov formula plays the key role in estimating the binary divisor sum \eqref{eq:dust}; the range of uniformity of $f$ with respect to $x$ depends on the best-known bound on the size of the Hecke eigenvalues of Maass forms (see \cite{Motohashi1}, Theorem 5).  Assuming the Ramanujan-Petersson conjecture would provide an asymptotic formula for \eqref{eq:divisorsumexample} in the region $N/M < q^{1 - \varepsilon}$.

Actually, our treatment is different than Motohashi's because we chose to use approximate functional equations in order to have formulas {\it a priori} valid in our region of interest (e.g., at the central point of the critical strip).  This feature causes the main terms to be captured in a more straightforward way in our work.  In addition, approximate functional equations explicitly display how large variables are with respect to each other.  Nevertheless, the basic strategy closely follows Motohashi.

One curious aspect of our work is an approximate generalization of Motohashi's formula for the fourth moment of the zeta function in terms of the third moment of Hecke-Maass L-functions.  We arrive at a certain formula involving products of three Hecke-Maass L-functions twisted by the $q$-th Hecke eigenvalue; see \eqref{eq:thirdmomentformula} below.  An important difference between these formulas is that our \eqref{eq:thirdmomentformula} is only helpful when $M$ and $N$ are restricted to be fairly close.


We have chosen to compute the shifted fourth moment of Dirichlet L-functions, that is we include the parameters $\alpha, \beta, \gamma, \delta$ similarly to \eqref{eq:motohashishift}.  Doing so allows for a more structural viewpoint of the main terms, and also allows us to check the conjecture of the five authors.

\subsection{Acknowledgements}
I am most grateful to Brian Conrey for extensive discussions on this material and for suggesting this problem.  I benefitted from conversations with David Farmer, Peng Gao, and Soundararajan.  I especially want to thank Henryk Iwaniec for a number of useful discussions.

\section{Notation and background}
\subsection{Dirichlet L-functions}
In this section we briefly recall some standard facts needed about Dirichlet L-functions.
Let $q$ be an odd prime ($q > 3$), $\chi$ be a primitive character modulo $q$, and $L(s,\chi)$ be the Dirichlet L-function
$L(s,\chi) = \sum_n \chi(n) n^{-s}$.  
The completed L-function satisfies the functional equation
\begin{equation*}
\Lambda(s, \chi) = \left(\frac{q}{\pi}\right)^{\frac{s}{2}} \Gamma\left(\frac{s + \sign}{2} \right) L(s, \chi)
 = 
\epsilon(\chi) \Lambda(1-s, \chibar),
\end{equation*}
where $\tau(\chi) = \sum_{x \pmod{q}} \chi(x) e\left(\frac{x}{q}\right)$ is the Gauss sum, and
\begin{equation*}
\epsilon(\chi) = i^{-\sign} q^{-\half} \tau(\chi),
\qquad
\sign = 
\begin{cases}
0, & \chi(-1) = 1 \\
1, & \chi(-1) = -1.
\end{cases}
\end{equation*}
In its asymmetric form, the functional equation reads
\begin{equation*}
L(s, \chi) = \epsilon(\chi) X(s) L(1-s, \chibar),
\end{equation*}
where
\begin{equation}
\label{eq:Xfactor}
X(\tfrac12 + u) = \left(\frac{q}{\pi} \right)^{-u } \frac{\Gamma\left(\frac{\half - u + \sign}{2} \right)}{\Gamma\left(\frac{\half + u + \sign}{2} \right)}.
\end{equation}

The central quantity of interest in this paper is
\begin{equation}
\label{eq:momentdefinition}
M(\alpha, \beta, \gamma, \delta) = \frac{2}{\phi^*(q)} \sideset{}{^+}\sum_{\chi \shortmod{q}} L(\tfrac12 + \alpha, \chi) L(\tfrac12 + \beta, \chi) L(\tfrac12 + \gamma, \chibar), L(\tfrac12 + \delta, \chibar),
\end{equation}
where the $+$ indicates the summation is restricted to primitive even characters and $\phi^*(q)$ is the number of primitive (odd or even) characters.  Throughout this work we assume that the shifts $\alpha, \beta, \gamma, \delta$ are all sufficiently small with respect to $q$ (say, $\ll(\log{q})^{-1}$).  It is natural to split the family separately into even characters and odd characters because the two families have different gamma factors in their functional equations.  In this work we concentrate almost exclusively on the even characters because the case of the odd characters is similar (we could treat both cases simultaneously but it would clutter the notation).  In Section \ref{section:odd} we briefly describe the necessary changes to treat the odd characters.

\subsection{The conjecture for the fourth moment}
\begin{myconj}{\cite{CFKRS}}
\label{conjecture:five}
For any $q \not \equiv 2 \pmod{4}$, and with shifts $\ll (\log{q})^{-1}$, 
\begin{multline}
\label{eq:conjecture}
M(\alpha, \beta, \gamma, \delta) 
= 
\frac{\zeta_q(1 + \alpha + \gamma)\zeta_q(1 + \alpha + \delta) \zeta_q(1 + \beta + \gamma)\zeta_q(1 + \beta + \delta)}{\zeta_q(2 + \alpha + \beta + \gamma + \delta)} \\
+ X_{\alpha, \gamma} \frac{\zeta_q(1 - \alpha +\beta)\zeta_q(1 - \alpha - \gamma) \zeta_q(1 + \beta + \delta)\zeta_q(1 - \gamma + \delta)}{\zeta_q(2 - \alpha + \beta - \gamma + \delta)} \\
+ X_{\alpha, \delta} \frac{\zeta_q(1 - \alpha + \beta)\zeta_q(1 - \alpha - \delta) \zeta_q(1 + \beta + \gamma)\zeta_q(1 + \gamma - \delta)}{\zeta_q(2 - \alpha + \beta + \gamma - \delta)} \\
+ X_{\beta, \gamma} \frac{\zeta_q(1 + \alpha -\beta) \zeta_q(1 + \alpha + \delta)\zeta_q(1 - \beta - \gamma) \zeta_q(1 - \gamma + \delta)}{\zeta_q(2 + \alpha - \beta - \gamma + \delta)} \\
+ X_{\beta, \delta} \frac{\zeta_q(1 + \alpha -\beta) \zeta_q(1 + \alpha + \gamma)\zeta_q(1 - \beta - \delta) \zeta_q(1 + \gamma - \delta)}{\zeta_q(2 + \alpha - \beta + \gamma - \delta)} \\
+ X_{\alpha, \beta, \gamma, \delta} \frac{\zeta_q(1 - \alpha - \gamma)\zeta_q(1 - \alpha - \delta) \zeta_q(1 - \beta - \gamma)\zeta_q(1 - \beta - \delta)}{\zeta_q(2 - \alpha - \beta - \gamma - \delta)} 
+ O(q^{-1/2 + \varepsilon}),
\end{multline}
where $\zeta_q(s) = \zeta(s) \prod_{p | q} (1 - p^{-s})$ and $X_{\alpha, \gamma} = X(\half + \alpha) X(\half + \gamma)$, etc.
\end{myconj}
In the case of $q$ prime one may replace each occurence of $\zeta_q$ with $\zeta$ without altering the error term.  There are no primitive characters modulo $q$ if $q \equiv 2 \pmod{4}$. 

The main terms on the right hand side clearly exhibit the same symmetries as the moment must (arising from trivial permutations of the variables as well as applications of the functional equation of the Dirichlet L-functions).  A similar conjecture holds for the odd characters; the only change is to have the $X$ factors given by \eqref{eq:Xfactor} depend on the parity.

It is important not to underestimate the psychological advantage of having the main terms produced ahead of time.  The actual computation of these main terms is by no means straightforward; there are many cancellations and combinations of various terms that occur throughout the work.  Although it may appear that \eqref{eq:conjecture} is a messy expression, in fact it arises after many simplifications.

\begin{mytheo}
 \label{thm:mainresultshifts}
For $q > 3$ prime, Conjecture \ref{conjecture:five} holds but with an error of size $O(q^{-5/512 + \varepsilon})$.
\end{mytheo}
Theorem \ref{thm:mainresult} follows by taking the limit as the shifts go to zero.  It is not obvious at first glance that the main term on the right hand side of \eqref{eq:conjecture} is holomorphic with respect to the shift parameters, but this is proven in a more general setting in \cite{CFKRS}; a computer can also easily verify it in this case.  Of course the left hand side is holomorphic in the shifts.

In the course of the work we assume that each of the shifts lies in a fixed annulus with inner and outer radii $\asymp (\log{q})^{-1}$; we further suppose that these annuli are separated enough so that any linear combination of the form $\epsilon_1 \alpha + \epsilon_2 \beta + \epsilon_3 \gamma + \epsilon_4 \delta$ with each $\epsilon_i \in \{ -1, 0, 1\}$, not all zero, is $\gg (\log{q})^{-1}$.  We initially prove Theorem \ref{thm:mainresultshifts} with this restriction in place; by the holomorphy of the left hand side of \eqref{eq:conjecture}, as well as the combination of main terms on the right hand side of \eqref{eq:conjecture}, we see that the error term also must be holomorphic in terms of the shifts.  By the maximum modulus principle, we can then extend the result to all shifts $\ll (\log{q})^{-1}$ with the same error term.

\subsection{The orthogonality formula}
\begin{myprop} 
\label{prop:orthogonality}
If $(ab, q) =1$ then
\begin{equation}
\sideset{}{^+}\sum_{\chi \shortmod{q}} \chi(a)\chibar(b) = \half \sum_{\substack{d | (q, a \pm b)}} \phi(d) \mu(q/d).
\end{equation}
The sum on the left hand side vanishes if $(ab, q) \neq 1$.  The condition $d | a \pm b$ should be taken with multiplicity (i.e., if $d|(a+b, a-b)$, it is counted twice).
\end{myprop}
The proof is standard and appears in many sources such as \cite{H-B1}, \cite{S}.  The odd characters have a similar orthogonality relation except the terms with $d| a + b$ are subtracted from the terms with $d | a-b$.  Also, note that for $q$ prime, $\phi^*(q) = q-2$, and the number of even, primitive characters is $\half(q-3)$.

\subsection{Approximate functional equation}
To get a useful formula for $M(\alpha,\beta,\gamma,\delta)$ we use an approximate functional equation.  There are a variety of choices in how to represent $|L(\half, \chi)|^4$; we have chosen to use the functional equation of $L(s,\chi)^2 L(s,\overline{\chi})^2$ to write $|L(\half, \chi)|^4$ as the sum of two sums of length approximately $q^2$.  Actually, we need a formula for $L(\half +\alpha, \chi) \dots L(\half + \delta, \overline{\chi})$, which prevents us from easily lifting a formula from, say, Theorem 5.3 of \cite{IK}, however the derivation of the formula follows standard lines.  We have
\begin{myprop}[Approximate functional equation]
\label{prop:AFE}
Let $G(s)$ be an even, entire function of exponential decay as 
$|s| \rightarrow \infty$ in any fixed strip $|\text{Re}(s)| \leq C$ and let
\begin{equation}
\label{eq:V}
V_{\alpha, \beta, \gamma, \delta}(x) = \frac{1}{2 \pi i} \int_{(1)} \frac{G(s)}{s} g_{\alpha, \beta, \gamma, \delta}(s) x^{-s} ds,
\end{equation}
where
\begin{equation}
g_{\alpha, \beta, \gamma, \delta}(s) = \pi^{-2s}
\frac{\Gamma\left(\frac{\half + \alpha + s + \sign}{2} \right)}{\Gamma\left(\frac{\half + \alpha  + \sign}{2} \right)} 
\frac{\Gamma\left(\frac{\half + \beta + s + \sign}{2} \right)}{\Gamma\left(\frac{\half + \beta  + \sign}{2} \right)} 
\frac{\Gamma\left(\frac{\half + \gamma + s + \sign}{2} \right)}{\Gamma\left(\frac{\half + \gamma  + \sign}{2} \right)}
\frac{\Gamma\left(\frac{\half + \delta + s + \sign}{2} \right)}{\Gamma\left(\frac{\half + \delta  + \sign}{2} \right)}.
\end{equation}
Furthermore, set
\begin{equation}
X_{\alpha,\beta,\gamma,\delta} = X(\tfrac12 + \alpha) X(\tfrac12 + \beta)X(\tfrac12 + \gamma)X(\tfrac12 + \delta).
\end{equation}
Then
\begin{multline}
L(\tfrac12 + \alpha, \chi) L(\tfrac12 + \beta, \chi) L(\tfrac12 + \gamma, \chibar)  L(\tfrac12+ \delta, \chibar) 
\\
=
 \sum_{m_1, m_2, n_1, n_2} \frac{\chi(m_1m_2) \chibar(n_1 n_2)}{m_1^{\half + \alpha}m_2^{\half + \beta} n_1^{\half + \gamma}n_2^{\half + \delta}} V_{\alpha, \beta, \gamma, \delta} \left(\frac{m_1m_2n_1n_2}{q^2}\right)
\\
+ X_{\alpha,\beta,\gamma,\delta}
   \sum_{m_1, m_2, n_1, n_2} \frac{\chibar(m_1m_2) \chi(n_1 n_2)}{m_1^{\half - \alpha}m_2^{\half - \beta} n_1^{\half - \gamma}n_2^{\half - \delta}} V_{-\alpha, -\beta, -\gamma, -\delta} \left(\frac{m_1m_2n_1n_2}{q^2}\right).
\end{multline}
\end{myprop}
We impose additional conditions on $G$ which we separate here for ease of reference.
\begin{mydefi}
\label{def:G}
Let $Q_{\alpha,\beta,\gamma,\delta}(s)$ be an even polynomial in $s$ with the following properties: it takes the value $1$ at $s=0$; it is rational in the shifts $\alpha, \beta, \gamma, \delta$; it is symmetric in the shifts; it is invariant under $\alpha \rightarrow -\alpha$, $\beta \rightarrow -\beta$, etc.; it is zero at $2s = -\alpha - \gamma$, $s = \half \pm \alpha$ (as well as $\beta + \delta$, $\half \pm \beta$, etc., by symmetry).  Then set $G(s) = Q_{\alpha,\beta,\gamma, \delta}(s) \exp(s^2)$.
\end{mydefi}

\begin{proof}[Proof of Proposition \ref{prop:AFE}]
Let
\begin{equation*}
\Lambda_{\alpha, \beta, \gamma, \delta}(s) = \Lambda(\tfrac12 + s + \alpha, \chi) \Lambda(\tfrac12 + s + \beta, \chi) \Lambda(\tfrac12 + s + \gamma, \chibar) \Lambda(\tfrac12 + s + \delta, \chibar)
\end{equation*}
and consider
\begin{equation*}
I_1 = \int_{(1)} \Lambda_{\alpha, \beta, \gamma, \delta}(s) \frac{G(s)}{s} ds.
\end{equation*}
Move the line of integration to $(-1)$, passing the pole at $s=0$.  Let $I_2$ be the new integral.  The residue at $s=0$ is
\begin{equation*}
\Lambda_{\alpha, \beta, \gamma, \delta}(0) = \Lambda(\tfrac12+ \alpha, \chi) \Lambda(\tfrac12 + \beta, \chi) \Lambda(\tfrac12 + \gamma, \chibar) \Lambda(\tfrac12 + \delta, \chibar).
\end{equation*}
After the change of variables $s \rightarrow -s$ and the application of the functional equation
\begin{equation*}
\Lambda_{\alpha, \beta, \gamma, \delta}(-s) = \Lambda_{-\gamma, -\delta, -\alpha, -\beta}(s),
\end{equation*}
we obtain
\begin{equation*}
I_2 = - \int_{(1)} \Lambda_{-\gamma, -\delta, -\alpha, -\beta}(s) \frac{G(s)}{s} ds.
\end{equation*}
Set 
\begin{equation*}
L_{\alpha, \beta, \gamma, \delta}(s) = L(\tfrac12 + \alpha +s, \chi) L(\tfrac12 + \beta +s, \chi) L(\tfrac12 + \gamma +s, \chibar) L(\tfrac12 + \delta +s, \chibar),
\end{equation*} 
and let
\begin{equation*}
\Lambda_{\alpha, \beta, \gamma, \delta}(s) = M_{\alpha, \beta, \gamma, \delta}(s) L_{\alpha, \beta, \gamma, \delta}(s).
\end{equation*}
Then we have
\begin{multline*}
L_{\alpha, \beta, \gamma, \delta}(0) = \frac{1}{2 \pi i} \int_{(1)} L_{\alpha, \beta, \gamma, \delta}(s) \frac{M_{\alpha, \beta, \gamma,\delta}(s)}{M_{\alpha, \beta, \gamma,\delta}(0)} \frac{G(s)}{s} ds 
\\
+ \frac{1}{2 \pi i} \int_{(1)} L_{-\gamma, -\delta, -\alpha, -\beta}(s) \frac{M_{-\gamma, -\delta, -\alpha, -\beta}(s)}{M_{\alpha, \beta, \gamma,\delta}(0)} \frac{G(s)}{s} ds.
\end{multline*}
An easy computation shows
\begin{equation*}
\frac{M_{\alpha, \beta, \gamma,\delta}(s)}{M_{\alpha, \beta, \gamma,\delta}(0)} = q^{2s} \left(\frac{q}{\pi}\right)^{2s}
  g_{\alpha, \beta, \gamma, \delta}(s),
  \quad
\frac{M_{-\gamma, -\delta, -\alpha, -\beta}(s)}{M_{\alpha, \beta, \gamma,\delta}(0)} 
= q^{2s} 
 X_{\alpha,\beta,\gamma,\delta} \,  g_{-\alpha, -\beta, -\gamma, -\delta}(s).
\end{equation*}
Expanding $L_{*, *, *, *}(s)$ into absolutely convergent Dirichlet series and reversing the order of summation and integration completes the proof.
\end{proof}

\subsection{Automorphic forms}
We briefly summarize the material we require on automorphic forms in order to apply the Kuznetsov formula.  See \cite{M}, \cite{Iwaniec2}, or \cite{IK} for further details.

Let $u_j(z)$ be an orthonormal system of Maass cusp forms on $SL_2(\mz) \backslash \mh$ with Laplace eigenvalues $\frac14 + \kappa_j^2$.  
Each cusp form $u_j(z)$ has the Fourier expansion
\begin{equation*}
u_j(z) = y^{\half} \sum_{n \neq 0} \rho_j(n) K_{i \kappa_j} (2 \pi |n| y) e(nx).
\end{equation*}
We assume that each $u_j$ is an eigenfunction of the Hecke 
operators $T_n$ defined by
\begin{equation*}
T_n u_j(z) = \frac{1}{\sqrt{n}} \sum_{ad=n} \sum_{b \shortmod{d}} u_j\left(\frac{az+b}{d}\right),
\end{equation*}
with eigenvalues $\lambda_j(n)$.  That is, $T_n u_j(z) = \lambda_j(n) u_j(z)$.  By consideration of the reflection operator, we may furthermore assume that $\rho_j(-n) = \epsilon_j \rho_j(n)$.  
The Hecke eigenvalues are multiplicative and satisfy the Hecke relation
\begin{equation}
\label{eq:heckerelation}
\lambda_j(m) \lambda_j(n) = \sum_{d|(m,n)} \lambda_j({\textstyle \frac{mn}{d^2} }),
\end{equation}
In terms of Fourier coefficients, the Hecke relations give $\rho_j(n) = \rho_j(1) \lambda_j(n)$.
Estimating the size of an individual Hecke eigenvalue is an important problem.  So far the best result is
\begin{equation}
|\lambda_j(n)| \leq d(n) n^{\theta},
\end{equation}
with $\theta = 7/64$, due to Kim and Sarnak \cite{KS}.  The Ramanujan-Petersson conjecture asserts that $\theta = 0$ is allowable.

The Hecke-Maass L-function is initially defined by
$L_j(s) = \sum_n \lambda_j(n) n^{-s}$,
which converges absolutely for $\text{Re}(s) > 1$ by properties of the Rankin-Selberg convolution.
Then $L_j(s)$ continues to an entire function and satisfies the functional equation
\begin{equation*}
\Lambda_j(s) := \pi^{-s} \Gamma \left(\frac{s + \delta_j + i\kappa_j}{2} \right)\Gamma \left(\frac{s + \delta_j - i\kappa_j}{2} \right) L_j(s) =\epsilon_j \Lambda_j(1-s),
\end{equation*}
where $\delta_j = (1 - \epsilon_j)/2$.
The Hecke relations \eqref{eq:heckerelation} translate to $L_j(s)$ having the 
Euler product
\begin{equation*}
L_j(s) = \prod_{p} \left( 1 - \frac{\lambda_j(p)}{p^s} + \frac{1}{p^{2s}} \right)^{-1}.
\end{equation*}

The Eisenstein series have the Fourier expansion
\begin{multline}
\label{eq:fouriereisenstein}
\pi^{-s} \Gamma(s) \zeta(2s) E(z,s) = \pi^{-s} \Gamma(s) \zeta(2s) y^s + \pi^{s-1} \Gamma(1-s) \zeta(2-2s) y^{1-s}
\\
+
2 y^{\half} \sum_{n \neq 0} |n|^{s-\half} \sigma_{1-2s}(|n|) K_{s-\half} (2 \pi |n| y) e(nx),
\end{multline}
where
\begin{equation}
\sigma_{\lambda} (n) = \sum_{d | n} d^{\lambda}.
\end{equation}
The Fourier expansion furnishes the meromorphic continuation of $E(z,s)$ to $s \in \mc$.

Let $u_{j,k}$, $1 \leq j \leq \dim{S_k}(SL_2(\mz))$ be a complete orthonormal Hecke basis of the classical weight $k$ cusp forms.  These have the Fourier expansion
\begin{equation*}
u_{j,k}(z) = \sum_{n \geq 1} \psi_{j,k}(n) n^{\frac{k-1}{2}} e(nz).
\end{equation*}
The Ramanujan-Petersson conjecture is known for holomorphic forms due to Deligne \cite{Deligne}.

\subsection{Kloosterman sums and Kuznetsov formula}
The Kuznetsov formula is a trace formula relating sums of Kloosterman sums to Fourier coefficients of automorphic forms.  This technology can show that there is considerable cancellation in the sum of Kloosterman sums.  Furthermore, it provides a separation of variables of $m$ and $n$ in $\sum_c S(m,n;c)$ which is conducive to obtaining additional savings in summations over $m$ and $n$; see \cite{IwaniecMean} and \cite{DI}.
\begin{mytheo}[Kuznetsov formula]
\label{thm:kuznetsov}
Let $g$ be a $C^2$ function satisfying $g(0) = 0$ and $g^{(j)}(x) \ll (x + 1)^{-2- \varepsilon}$, $j=0,1,2$ and suppose $m, n > 0$.  Then
\begin{multline}
\sum_{c}  \frac{S(m,n;c)}{c} g \left(\frac{4 \pi \sqrt{mn}}{c}\right)
= 
\sum_j \rho_{j}(m)\overline{\rho_{j}(n)} \mathcal{M}_g(\kappa_j)
\\
+ 
\sum_{k \equiv 0 \shortmod{2}} \mathcal{N}_g(k) \sum_j  \psi_{jk}(m) \overline{\psi}_{jk}(n)
+ \frac{1}{\pi}  \int_{-\infty}^{\infty} \cosh(\pi t) \mathcal{M}_g(t) 
 \frac{\sigma_{2it}(m) \sigma_{2it}(n)}{(mn)^{it} |\zeta(1 + 2it)|^2} dt.
\end{multline}
Here $\mathcal{M}_g$ and $\mathcal{N}_g$ are the following integral transforms
\begin{align*}
\mathcal{M}_g(t) &= \frac{\pi i}{\sinh{2\pi t}} \int_0^{\infty} (J_{2it}(x) - J_{-2it}(x)) g(x) \frac{dx}{x}, 
\\
\mathcal{N}_g(k) &= \frac{4 (k-1)! }{(4 \pi i)^{k}} \int_0^{\infty} J_{k-1}(x) g(x) \frac{dx}{x}.
\end{align*}
For the opposite sign case,
\begin{multline}
\sum_{c}  \frac{S(m,-n;c)}{c} g \left(\frac{4 \pi \sqrt{mn}}{c}\right)
=
\sum_j \rho_{j}(m)\overline{\rho_{j}(-n)} \mathcal{K}_{g}(\kappa_j)
\\
+ \frac{1}{\pi}  \int_{-\infty}^{\infty} \cosh(\pi t) \mathcal{K}_g(t) 
 \frac{\sigma_{2it}(m) \sigma_{2it}(n)}{(mn)^{it} |\zeta(1 + 2it)|^2} dt,
\end{multline}
where
\begin{equation*}
\mathcal{K}_{g}(t) = 2 \int_0^{\infty} K_{2it}(x) g(x) \frac{dx}{x}.
\end{equation*}
\end{mytheo}
For proofs we refer to Theorems 16.5 and 16.6 of \cite{IK} or Theorems 2.3 and 2.4 of \cite{M} (we borrowed some notation from both sources).

\subsection{The Estermann $D$-function}
The Estermann $D$-function is defined by
\begin{equation}
D(s, \lambda, \frac{h}{l} ) = \sum_{n} \frac{\sigma_{\lambda}(n)}{n^s} e\left(n\frac{h}{l}\right),
\end{equation}
where $(h,l) =1$.  The analytic properties of $D$ are useful for understanding the behavior of the divisor function.  We have
\begin{mylemma}
\label{lem:Dq}
For any fixed $\lambda \in \mc$, $D(s,\lambda, \frac{h}{l})$ is meromorphic as a function of $s$, and satisfies the functional equation
\begin{multline}
\label{eq:FE}
D(\tfrac12 +s, \lambda, \frac{h}{l}) = 2 (2\pi)^{-1 -\lambda +2s }  \Gamma(\tfrac12-s) \Gamma(\tfrac12+\lambda - s)
l^{\lambda - 2s}
\\
\times 
\Big[
D(\tfrac12-s, -\lambda, \frac{\overline{h}}{l} ) \cos(\tfrac{\pi \lambda}{2}) 
+ D(\tfrac12-s, -\lambda, - \frac{\overline{h}}{l}) \sin(\pi(s - \tfrac{\lambda}{2}))
\Big].
\end{multline}
If $\lambda \neq 0$ then $D$ has simple simple poles at $s=1$ and $s = 1 + \lambda$ with respective residues
\begin{equation}
 l^{-1 + \lambda}  \zeta(1 - \lambda), \qquad
 l^{-1 - \lambda}  \zeta(1 + \lambda).
\end{equation}
\end{mylemma}
We refer to \cite{M}, Lemma 3.7 for proofs.  The functional equation of the Estermann function is essentially equivalent to the Voronoi summation formula (Theorem 4.10 of \cite{IK}).

\subsection{Conventions}
We use the common convention in analytic number theory that $\varepsilon$ denotes an arbitrarily small positive quantity that may vary from line to line.  Furthermore, $\alpha$, $\beta$, $\gamma$, and $\delta$ are complex numbers that are sufficiently small in comparison to $\varepsilon$ (so for example we may say $\zeta(1 + \alpha + s)$ is holomorphic for $\text{Re}(s) > \varepsilon$).

\section{The structure of the fourth moment}
\subsection{Averaging the approximate functional equation}
Using Propositions \ref{prop:orthogonality} and \ref{prop:AFE} we may average the approximate functional equation to obtain a formula for the fourth moment $M(\alpha,\beta,\gamma,\delta)$ defined by \eqref{eq:momentdefinition}.  Write
\begin{equation}
\label{eq:MAFE}
M(\alpha, \beta, \gamma, \delta) = A_{1, q}(\alpha, \beta, \gamma, \delta) + A_{-1, q}(\alpha, \beta, \gamma, \delta),
\end{equation}
where $A_{1, q}$ is the contribution from the `first part' of the approximate functional equation, and $A_{-1, q}$ is the second part.  It suffices to compute $A_{1, q}$ since
\begin{equation}
\label{eq:switch}
A_{-1, q}(\alpha, \beta, \gamma, \delta) = X_{\alpha,\beta,\gamma,\delta} A_{1, q}(-\alpha, -\beta, -\gamma, -\delta).
\end{equation} 
We have
\begin{equation*}
A_{1, q} = \frac{1}{\phi^{*}(q) }\sum_{d | q} \phi(d) \mu(q/d) \mathop{\sum \sum \sum \sum}_{\substack{m_1 m_2 \equiv \pm n_1 n_2 \shortmod{d} \\ (m_1 m_2 n_1 n_2, q) = 1} } \frac{1}{m_1^{\half + \alpha} m_2^{\half + \beta} n_1^{\half + \gamma} n_2^{\half + \delta}} V_{\alpha, \beta, \gamma, \delta}\left(\frac{m_1m_2n_1n_2}{q^2}\right).
\end{equation*}
By changing variables $m = m_1 m_2$ and $n = n_1 n_2$, we may rewrite this as
\begin{equation}
\label{eq:A1useful}
A_{1, q} = \frac{1}{\phi^{*}(q)} \sum_{d | q} \phi(d) \mu(q/d) \mathop{\sum \sum}_{\substack{m \equiv \pm n \shortmod{d} \\ (m n, q) = 1} } \frac{\sigma_{\alpha - \beta}(m) \sigma_{\gamma - \delta}(n)}{m^{\half + \alpha} n^{\half + \gamma} } V_{\alpha, \beta, \gamma, \delta}\left(\frac{mn}{q^2}\right).
\end{equation}
Let $A_1$ be the same sum as $A_{1, q}$ but with the condition $(mn,q) = 1$ omitted.  Since $q$ is prime, a trivial estimation gives
\begin{equation}
\label{eq:A1q}
A_{1, q} = A_{1} + O(q^{-1 + \varepsilon}).
\end{equation}

\subsection{The diagonal terms}
In this section we compute the diagonal contribution $A_D$ of the terms $m=n$ in $A_1$ given by \eqref{eq:A1useful}.  To be precise we are computing the diagonal contribution of the terms where $m \equiv n \pmod{d}$ and not those with $m \equiv - n \pmod{d}$.  We compute
\begin{align*}
A_D &= \frac{1}{\phi^{*}(q)}  \sum_{d|q} \phi(d) \mu(q/d) \sum_{(n, q) = 1} \frac{\sigma_{\alpha - \beta}(n) \sigma_{\gamma - \delta}(n)}{ n^{1 + \alpha + \gamma} } V_{\alpha, \beta, \gamma, \delta}\left(\frac{n^2}{q^2}\right) 
\\
&= \sum_{(n, q) = 1} \frac{\sigma_{\alpha - \beta}(n) \sigma_{\gamma - \delta}(n)}{ n^{1 + \alpha + \gamma} } V_{\alpha, \beta, \gamma, \delta}\left(\frac{n^2}{q^2}\right),
\end{align*}

Now the problem reduces to a standard exercise in analytic number theory.  We use the integral representation of $V$, that is \eqref{eq:V}, to write the sum as an integral of a Dirichlet series and develop the asymptotics by moving the line of integration.  Precisely, we have
\begin{equation*}
A_D =  \frac{1}{2 \pi i} \int_{(1)} \frac{G(s)}{s} g_{\alpha, \beta, \gamma, \delta}(s) q^{2s} 
 \sum_{(n, q) = 1} \frac{\sigma_{\alpha - \beta}(n) \sigma_{\gamma - \delta}(n)}{ n^{1 + \alpha + \gamma + 2s} } ds.
\end{equation*}
From the Ramanujan identity
\begin{equation}
\label{eq:Ramanujan}
\sum_n \frac{\sigma_{\lambda}(n) \sigma_{\nu}(n)}{n^v} = \frac{\zeta(v) \zeta(v - \lambda) \zeta(v - \nu) \zeta(v - \lambda - \nu)}{\zeta(2v - \lambda - \nu)}
\end{equation}
we obtain
\begin{equation*}
A_D =  \frac{1}{2 \pi i} \int_{(1)} \frac{G(s)}{s} g_{\alpha, \beta, \gamma, \delta}(s) q^{2s} Z_{\alpha, \beta, \gamma, \delta}(s) ds,
\end{equation*}
where
\begin{equation*}
Z_{\alpha, \beta, \gamma, \delta}(s) = \frac{\zeta(1 + \alpha + \gamma + 2s) \zeta(1 + \alpha + \delta + 2s)\zeta(1 + \beta + \gamma + 2s) \zeta(1 + \beta + \delta + 2s) }{\zeta(2 + \alpha + \beta + \gamma + \delta + 4s)}.
\end{equation*}
$Z$ has simple poles at $2s = -\alpha - \gamma$, $-\alpha - \delta$, $-\beta - \gamma$, and $-\beta - \delta$ and is otherwise holomorphic for $\text{Re } s > -1/4 + \varepsilon$.  Note that $G(s)$ vanishes at these poles of $Z$ by our choice in Definition \ref{def:G}.  We move the line of integration to the line $-1/4 + \varepsilon$, passing the pole at $s=0$ only.  We obtain $A_D = (\text{Residue}) + I$, where $I$ is the contribution from the new line of integration.  
We easily obtain
\begin{equation*}
|I| \ll q^{-\half + \varepsilon}.
\end{equation*}
The residue at $s=0$ gives
\begin{equation}
Y_1(\alpha, \beta; \gamma, \delta) :=  \frac{\zeta(1 + \alpha + \gamma) \zeta(1 + \alpha + \delta)\zeta(1 + \beta + \gamma ) \zeta(1 + \beta + \delta ) }{\zeta(2 + \alpha + \beta + \gamma + \delta )},
\end{equation}
which is precisely the first term in Conjecture \ref{conjecture:five} (up to $O(q^{-1 + \varepsilon})$).  
We summarize this calculation with the following
\begin{mylemma}
\label{lem:YandJ}
With $G$ as in Definition \ref{def:G}, we have
\begin{equation}
\label{eq:YandJ}
A_D = Y_1(\alpha, \beta; \gamma, \delta) 
+ 
O(q^{-1/2 + \varepsilon}),
\end{equation}
and similarly the contribution of the diagonal terms to $A_{-1}$ is
\begin{equation}
\label{eq:YandJ2}
A_{-D} = Y_{-1}(\alpha, \beta; \gamma, \delta) 
+ O(q^{-1/2 + \varepsilon}).
\end{equation}
Here
\begin{align}
Y_{-1}(\alpha, \beta; \gamma, \delta) &= X_{\alpha,\beta,\gamma,\delta} Y_1(-\alpha, -\beta;-\gamma, -\delta).
\end{align}
\end{mylemma}
If we did not choose $G$ to vanish at the poles of $Z$ then there would be four more extra terms in each of \eqref{eq:YandJ} and \eqref{eq:YandJ2}.  However, these extra terms would all cancel out in the final assembly of $M(\alpha,\beta,\gamma,\delta)$; this satisfying (but complicated) calculation is carried out in the first version of this article on the arxiv.

\subsection{Decomposition of the off-diagonal terms}
Here we investigate the contribution $A_{O}$ to $A_1$ from the off-diagonal terms $m \equiv n \pmod{d}$, $m \neq n$.  The treatment of the `dual' terms $A_{\overline{O}}$ corresponding to $m \equiv -n \pmod{d}$ proceeds in much the same way but must be executed separately, and is carried out in Section \ref{section:dual}.

We need to compute
\begin{equation}
A_O = \frac{1}{\phi^{*}(q)} \sum_{d | q} \phi(d) \mu(\frac{q}{d}) \mathop{\sum \sum {}'}_{\substack{m \equiv n \shortmod{d}} } \frac{\sigma_{\alpha - \beta}(m) \sigma_{\gamma - \delta}(n)}{m^{\half + \alpha} n^{\half + \gamma} } V_{\alpha, \beta, \gamma, \delta}\left(\frac{mn}{q^2}\right),
\end{equation}
where the prime indicates that the case $m=n$ is excluded.  

We break the sum $A_O$ into two pieces depending on whether $m < n$ or $m > n$ and write $A_O = B_O + B_O'$ accordingly.  By symmetry, it suffices to compute $B_O = B_O(\alpha, \beta, \gamma, \delta)$, since $B_O'(\alpha, \beta, \gamma, \delta) = B_O(\gamma, \delta, \alpha, \beta)$.  We record the decomposition,
\begin{equation}
\label{eq:decomposition}
 A_1 = A_D + B_O + B_O' + A_{\overline{O}}.
\end{equation}
where each quantity has parameters $\alpha, \beta, \gamma, \delta$ (in that order).

We require a partition of unity with some special properties.  To begin, we take a partition of unity $\{ W_M(x) \}$ on $\mr^+$ where each $W_M(x)$ has support in the dyadic interval $[M, 2M]$, and furthermore $W_M(x) = W(\frac{x}{M})$ where $W$ is a fixed smooth, compactly supported function.  Here $M$ varies over a set of positive real numbers, with the number of such $M$ less than $X$ being $O(\log{X})$.  Such a partition can be constructed very easily.  We then create the partition of unity $\{W_{M,N}(x,y) \}$ on $\mr^+ \times \mr^+$ by taking products of the $W_M(x)$.

We apply this partition of unity to $B_{O}$ and write  $B_{O} = \sum_{M, N} B_{M,N}$, where $B_{M,N}$ is the same expression as $B_O$ but weighted with $W_{M,N}(m,n)$.  Since $m < n$ we may assume 
\begin{equation}
M \ll N,
\end{equation}
a convention that holds throughout this paper.  Due to the rapid decay of $V(x)$, we may assume $MN \leq q^{2 + \varepsilon}$.

\subsection{Overview of the proof}
We write
\begin{equation}
B_{M,N} = (\text{Main term})_{M,N} + E_{M,N},
\end{equation}
for a certain main term that we do not explicitly write here due to its complexity.  The discussion of the main terms takes place primarily in Sections \ref{section:mainterms} and \ref{section:assemble}.  In Section \ref{section:upperboundonP} we prove that the size of the main term as a function of $M$ and $N$ is given by
\begin{equation}
\label{eq:maintermbound}
(\text{Main term})_{M,N} \ll M^{\half} N^{-\half} q^{\varepsilon},
\end{equation}
so it is a bit of a misnomer to call this a `main term' when $M$ and $N$ are not near each other.

In Section \ref{section:divisorAP} we prove the following estimate for $B_{M,N}$ that is applicable when $N$ is significantly larger than $M$.
\begin{mytheo}
\label{prop:faraway}
For $M \ll N$, $MN \ll q^{2+\varepsilon}$, we have
\begin{equation}
\label{eq:faraway}
E_{M,N} \ll 
(M^{\half} N^{-\frac14} + q^{\frac14} N^{-\frac14} + q^{\frac{3}{10}} N^{-\frac16} M^{-\frac{2}{15}}) q^{\varepsilon}.
\end{equation}
\end{mytheo}
Actually we prove the bound \eqref{eq:faraway} but with $B_{M,N}$ replacing $E_{M,N}$.
Notice the first term above is larger than the bound \eqref{eq:maintermbound} for the `main term', so in fact the two bounds are equivalent.  We do not attempt to extract a main term in this analysis, but still safely claim the same main term relevant in the range where $M$ and $N$ are close.

To cover the range where $M$ and $N$ are fairly close, we prove the following in Section \ref{section:errorterms}.
\begin{mytheo} 
\label{thm:Ebound}
With $\theta = 7/64$, $M \ll N$, and $MN \ll q^{2+\varepsilon}$, we have
\begin{equation}
\label{eq:Ebound}
E_{M,N} \ll q^{-\half + \theta + \varepsilon} M^{-\half} N^{\half}.
\end{equation}
\end{mytheo}

The various main terms combine in a rather complicated way.  When summed over all $M$ and $N$, all the main terms from $A_O$, $A_{\overline{O}}$, etc. add up to form the quantity in Conjecture \ref{conjecture:five}.

For some small values of $MN$ we can do no better than the trivial bound
\begin{equation}
\label{eq:trivial}
B_{M,N} \ll \frac{\sqrt{MN}}{q} q^{\varepsilon}.
\end{equation}
Note again that this bound on $B_{M,N}$ is larger than the bound \eqref{eq:maintermbound} for the main term provided $N > q$ (which is trivially satisfied since $B_{M,N}$ is void otherwise), so this translates into a bound on $E_{M,N}$.

Deducing Theorem \ref{thm:mainresultshifts} from the above results
is then an exercise in finding the maximum of a piecewise linear function.  Precisely, let
\begin{gather*}
L_1(a,b) = -1 + \half a + \half b, \qquad
L_2(a,b) = -\half + \theta -\half a + \half b, \\
L_3(a,b) = \half a - \frac14 b, \qquad
L_4(a,b) = \frac14 - \frac14 b, \qquad
L_5(a,b) = \frac{3}{10} - \frac{2}{15} a - \frac16 b.
\end{gather*}
We require
\begin{equation*}
\max_{i \in \{3,4,5\} } \max_{\substack{0 \leq a \leq b \\ a+ b \leq 2}} \min(L_1(a,b), L_2(a,b), L_i(a,b)).
\end{equation*}
Along the line $a=0$ the maximum value is $-\frac{1}{40}$ at $b=\frac{39}{20}$.  When $a=b$ the maximum is $-\half + \theta$ and along the line $a+ b= 2$ the maximum is $-\frac{1}{62} + \frac{\theta}{31}$ (although if $\theta$ was not known to be less than $6/43$ then $-\frac{1}{14} + \frac{3\theta}{7}$ would be the maximum here).  The maximum at an interior point must occur when $L_1 = L_2 = L_i$ for $i=3$, $4$, or $5$.  The case $i=3$ is $a = \half + \theta$, $b = \frac{4}{3}$, which gives the bound $-\frac{1}{12} + \frac{\theta}{2}$.  The case $i=4$ gives $-\frac{1}{12} + \frac{\theta}{6}$ from $a = \half + \theta$, $b = \frac43 - \frac{2}{3} \theta$ .  The case $i=5$ has $a=\half + \theta$, $b=\frac{59}{40} - \frac{19}{20} \theta$, and gives $-\frac{1}{80} + \frac{\theta}{40}$, completing the proof.  Note that our method requires the bound $\theta < 1/6$ which is very deep.

\section{The divisor function in arithmetic progressions}
\label{section:divisorAP}
In this section we present a handful of different estimates that taken together prove Theorem \ref{prop:faraway}.  We largely use classical techniques of analytic number theory such as Poisson summation, Cauchy's inequality, `gluing' of variables to create a longer variable, reciprocity laws, the Weil bound for Kloosterman sums, estimates for exponential sums, etc.

Solely for notational convenience we shall work with the case $\alpha = \beta = \gamma = \delta = 0$; the arguments extend easily to nonzero parameter values.

\subsection{Initial cleaning and statement of results}
Recall
\begin{equation}
\label{eq:lipton}
B_{M,N} = \frac{1}{\phi^{*}(q)} \sum_{d | q} \phi(d) \mu(q/d) \mathop{\sum \sum {}'}_{\substack{m \equiv n \shortmod{d}} } \frac{d(m) d(n)}{m^{\half} n^{\half} } V\left(\frac{mn}{q^2}\right) W\left(\frac{m}{M}\right) W\left(\frac{n}{N}\right).
\end{equation}
Throughout this section, $W(x)$ will denote a smooth function with support in a dyadic interval, which may vary from line to line (simply to avoid cluttering the notation).  Since we treat the range $M \asymp N$ in a different way, we assume $N > cM$ for a sufficiently large constant $c$ so that the term $m=n$ is avoided by the support of the weight functions in \eqref{eq:lipton}.
As a simple first step, we separate the variables $m$ and $n$ in the test function $V$ by using its Mellin transform.  We have
\begin{multline*}
B_{M,N} = \frac{1}{2 \pi i} \int_{(\varepsilon)} \left(\frac{q^2}{MN}\right)^{s} \frac{G(s)}{s}g(s)
\\
\frac{1}{\phi^{*}(q)} \sum_{d | q} \phi(d) \mu(q/d) \mathop{\sum \sum}_{\substack{m \equiv n \shortmod{d}} } \frac{d(m) d(n)}{m^{\half} n^{\half} }  W_s\left(\frac{m}{M}\right) W_s\left(\frac{n}{N}\right)  ds,
\end{multline*}
where $W_s(x) = x^{-s} W(x)$.  Note that the dependence of $W_s$ on $s$ is very mild:
\begin{equation*}
x^{j} \frac{d^j}{dx^j} W_s(x) \ll x^{-\text{Re}(s)} |P_j(s)| \max_{0 \leq i \leq j} \Big|x^i \frac{d^i}{dx^i} W(x)\Big|,
\end{equation*}
where $P_j(s)$ is a degree $j$ polynomial in $s$.  In effect, this separation of variables has no cost since losses of size $(MN)^{\varepsilon}$ are absorbed by the factor $q^{\varepsilon}$ in the bound \eqref{eq:faraway}, and because $s^{-1} G(s) g(s)$ has exponential decay in the imaginary direction.  Thus, the problem reduces to bounding $B'_{M,N}$, which is defined to be the same sum as $B_{M,N}$ but with $V$ removed.  Let
\begin{equation*}
\Delta_{M,N}(m) = \frac{1}{\phi^{*}(q)} \sum_{d | q} \phi(d) \mu(q/d) \sum_{n \equiv m \shortmod{d}}  d(n)  W\left(\frac{n}{N}\right),
\end{equation*}
so that
\begin{equation*}
B'_{M,N} = \frac{1}{\sqrt{MN}} \sum_{m} d(m) \Delta_{M,N}(m)  W\left(\frac{m}{M}\right),
\end{equation*}
where we have redefined $W(x)$ again to include the scaling factor $x^{-\half}$.

To begin, we open the divisor function $d(n)$ and let $n = n_1 n_2$.  To locate the variables $n_1$ and $n_2$, we apply
 a dyadic partition of unity to both $n_1$ and $n_2$ so that $n_1 \asymp N_1$ and $n_2 \asymp N_2$, where $N_1 N_2 \asymp N$. Without loss of generality assume $N_2 \geq N_1$.  Let $B_{M, N_1,N_2}$ be the sum $B'_{M,N}$ but with this partition of unity applied, and similarly for $\Delta_{M,N_1,N_2}(m)$.  

We prove the following in Section \ref{section:cauchy}
\begin{mylemma}
\label{lemma:redplastic}
Suppose $N_2 \geq N_1$, $N_1 N_2 \asymp N \gg M$, $MN \ll q^{2+\varepsilon}$.  Then
\begin{equation}
\label{eq:redplastic}
B_{M,N_1,N_2} \ll 
\begin{cases}
\left(\frac{N}{M}\right)^{-\half} q^{\half + \varepsilon} \\
(q^{\frac14} N_1^{-\half}  + \left(\frac{N}{M}\right)^{-\half} N_1^{\half}) q^{\varepsilon} \\
\left(\frac{N}{M}\right)^{-\half} N_1 q^{\varepsilon}.
\end{cases}
\end{equation}
\end{mylemma}
The first bound is nontrivial for $ N M^{-1} > q^{1 + \varepsilon}$ while the second bound is nontrivial for $q^{\half + \varepsilon} < N_1 < N M^{-1} q^{-\varepsilon}$ (note this upper bound is achieved automatically if $N \gg M^2 q^{\varepsilon}$ since $N_1 \ll \sqrt{N}$).  The third bound is nontrivial for $N_1 < N^{\half} M^{-\half}q^{-\varepsilon}$.  Note that if $N M^{-1} < q$, then there is a gap that has not been covered.  In Section \ref{section:N1small} we treat the case where $N_1$ is relatively small, proving
\begin{mylemma}
\label{lemma:chewyquaker}
Suppose $N_2 \geq N_1$, $N_1 N_2 \asymp N \gg M$, $MN \ll q^{2+\varepsilon}$. Then
\begin{equation}
\label{eq:chewyquaker}
B_{M,N_1,N_2} \ll \left(N_1^2 \min\left(\frac{q^{\half}}{N^{\frac56} M^{\frac23}}, \frac{1}{\sqrt{MN}}\right) + M^{\half} N^{-\half}  + M^{\frac32} N^{-\frac32} N_1 \right)q^{\varepsilon}.
\end{equation}
\end{mylemma}
We briefly describe how to deduce Theorem \ref{prop:faraway} from Lemmas \ref{lemma:redplastic} and \ref{lemma:chewyquaker}.  Basically we need to find the maximum value of the minimum of the above bounds as $N_1$ varies between $1$ and $\sqrt{N}$.  The case $N_1 = 1$ gives $\ll N^{-\half} M^{\half} q^{\varepsilon}$ and $N_1 = \sqrt{N}$ gives $\ll (M^{\half} N^{-\frac14} + q^{\frac14} N^{-\frac14}) q^{\varepsilon}$ (which already dominates the $N_1 = 1$ contribution).  Note that $M^{\half} N^{-\half} + M^{\frac32} N^{-\frac32} N_1$ is dominated by $N^{\frac14} M^{-\half}$ provided $M < N^{\half}$, which we may assume since otherwise \eqref{eq:faraway} is trivial. Hence the problem reduces to finding
\begin{equation*}
\max_{1 \leq N_1 \leq \sqrt{N}}(\min(c_1 N_1^{-\half},  c_2 N_1^2)), \quad \text{where }
c_1 = q^{\frac14}, 
\quad c_2 = \min\left(\frac{q^{\half}}{N^{\frac56} M^{\frac23}}, \frac{1}{\sqrt{MN}}\right).
\end{equation*}
The maximum is clearly
\begin{equation*}
q^{\frac15} \min\left(\frac{q^{\half}}{N^{\frac56} M^{\frac23}}, \frac{1}{\sqrt{MN}}\right)^{\frac15}
\ll
\frac{q^{\frac{3}{10}}}{N^{\frac16} M^{\frac{2}{15}}}.
\end{equation*}

\subsection{Proof of Lemma \ref{lemma:redplastic}}
\label{section:cauchy}
We begin by noting that we may assume $(n,d) = 1$ since the contribution to $B_{M,N}$ from $q|n$ is $O(q^{-1+\varepsilon})$ which is smaller than all the bounds in \eqref{eq:redplastic} as well as \eqref{eq:chewyquaker}.  Define $\Delta'$ to be the same as $\Delta$ but with $(n,d) = 1$.

Next we break up the sum over $n_2$ into arithmetic progressions $\pmod{d}$, getting
\begin{equation}
\label{eq:delta}
\Delta'_{M,N_1,N_2}(m) = \frac{1}{\phi^{*}(q)} \sum_{d | q} \phi(d) \mu(q/d) \sum_{(n_1, d)=1} \sum_{n_2 \equiv m \overline{n_1} \shortmod{d}} W\Big(\frac{n_1 n_2}{N}\Big) W\Big(\frac{n_1}{N_1}\Big) W\Big(\frac{n_2}{N_2}\Big) 
\end{equation}
Using the same argument we used to remove $V$ reduces the problem to estimating $\Delta''_{M,N_1,N_2}(m)$, where $\Delta''$ has the same form as $\Delta'$ but with $W(n_1n_2/N)$ removed.  Applying Poisson summation to the summation over $n_2$ gives
\begin{equation*}
\Delta''_{M,N_1,N_2}(m) = \frac{1}{\phi^{*}(q)} \sum_{d | q} \phi(d) \mu(q/d) \sum_{(n_1,d)=1} \frac{N_2}{d} \sum_{h} e\left(\frac{hm \overline{n_1}}{d}\right) W\left(\frac{n_1}{N_1}\right) \widehat{W}\left(\frac{h}{d/N_2}\right).
\end{equation*}
where $\widehat{W}$ is the Fourier transform of $W$.
The term $h=0$ contributes to $\Delta''_{M,N_1,N_2}$ precisely
\begin{equation*}
\frac{N_2}{\phi^{*}(q)} \widehat{W}(0) \sumstar_{n_1}  W\left(\frac{n_1}{N_1}\right) 
\sum_{d | q} \frac{\phi(d)}{d} \mu(q/d) \ll \frac{N}{q^2} q^{\varepsilon},
\end{equation*}
and hence its contribution to $B'_{M,N}$ is
\begin{equation*}
\ll \frac{\sqrt{MN}}{q^2} q^{\varepsilon} \ll q^{-1 + \varepsilon}.
\end{equation*}
From the rapid decay of $\widehat{W}$, the terms with $d=1$, $h \neq 0$ contribute $O(q^{-1000})$ to $B_{M,N_1, N_2}$.
It remains to consider the terms with $d=q$, $h \neq 0$, namely
\begin{equation}
\label{eq:brain}
R = \frac{N_2 \phi(q)}{q \phi^*(q) \sqrt{MN}}  \sum_{m} d(m) W\left(\frac{m}{M}\right) \sumstar_{n_1} \sum_{h \neq 0} e\left(\frac{hm \overline{n_1}}{q}\right)W\left(\frac{n_1}{N_1}\right) \widehat{W}\left(\frac{h}{H}\right),
\end{equation}
where $H = q N_2^{-1}$.
Due to rapid decay of the Fourier transform, we may assume $|h| \leq H q^{\varepsilon}$ (which implies $N_2 \leq q^{1+\varepsilon}$).  Furthermore we may assume $(hm, q) = 1$ at no cost.

It is possible to bound the sum over $n_1$ using the Weil bound for Kloosterman sums, however one can obtain better bounds on average over $m$ and $h$.  We glue together $m$ and $h$ to create a longer variable $l = mh$.  We have
\begin{equation}
\label{eq:greenmarker}
R \ll \frac{N_2 q^{\varepsilon}}{q \sqrt{MN}} 
\sumstar_{0<l \leq L}  \left| \sumstar_{n_1}  e\left(\frac{l \overline{n_1}}{q}\right) W\left(\frac{n_1}{N_1}\right) \right|,
\end{equation}
where $L \asymp MH$.  Bounding this double sum is of independent interest.  Friedlander and Iwaniec (\cite{FriedlanderI}, p. 337) state a bound for this sum which is nontrivial for $N_1$ and $L$ rather short (e.g., $L = q^{1/8}$ and $N_1 = q^{3/8 + \varepsilon}$), however we do not use their bound.
We have
\begin{myprop}
\label{thm:expsum}
Let $W$ be a smooth function with support in $[1,2]$, and let
\begin{equation}
S(K, L;q) := \sumstar_{0<l \leq L}  \left| \sumstar_{k}  e\left(\frac{l \overline{k}}{q}\right) W\left(\frac{k}{K}\right) \right|.
\end{equation}
If $K, L \ll q^{1 + \varepsilon}$, then
\begin{equation}
\label{eq:SKLq}
S(K, L;q) \ll
\begin{cases}
L q^{\half + \varepsilon}, \\
(L^{\half} q^{\frac{3}{4}} +  K^{\half} L) q^{\varepsilon}. 
\end{cases}
\end{equation}
\end{myprop}
Delaying the proof of Proposition \ref{thm:expsum} for a moment, we complete the proof of Lemma \ref{lemma:redplastic} by inserting \eqref{eq:SKLq} as well as the trivial bound into \eqref{eq:greenmarker}.

In the course of proving Proposition \ref{thm:expsum}, we need the following.
\begin{mylemma}
 \label{lemma:Tsum}
For any odd prime $q$, and $x, y, z, \in \mf_q$, let
\begin{equation}
T(x,y,z;q) := \sumstar_{a \shortmod{q}} \sumstar_{b \shortmod{q}} \sumstar_{c \shortmod{q}}
e\left(\frac{c(\overline{a} - \overline{b})}{q} \right)
e\left(\frac{a x + b y + c z}{q} \right).
\end{equation}
Then we have the following bounds,
\begin{equation}
\label{eq:Tbound}
T(x,y,z;q) \ll
\begin{cases}
q, \qquad z = 0, x \neq -y, \\
q^2, \qquad z = 0, x = -y, \\
q, \qquad z \neq 0, xy=0, \\
q^{\frac32}, \qquad x,y,z \neq 0.
\end{cases}
\end{equation}
\end{mylemma}
Actually the proof gives an exact formula for $T$ in terms of Kloosterman sums.
\begin{proof}[Proof of Lemma \ref{lemma:Tsum}]
 First note that the sum over $c$ is $-1$ unless $\overline{a} - \overline{b} + z \equiv 0$, in which case it equals $q-1$.  Hence
\begin{equation*}
T(x,y,z;q) = (q -1) \mathop{\sumstar \sumstar}_{\overline{a} - \overline{b} \equiv -z \shortmod{q}} e\left(\frac{a x + b y }{q} \right) - \mathop{\sumstar \sumstar}_{\overline{a} - \overline{b} \not \equiv -z \shortmod{q}} e\left(\frac{a x + b y }{q} \right).
\end{equation*}
Solving for $b$ and using $\overline{\overline{a} + z} = a (\overline{1 + z a})$ gives
\begin{equation*}
T(x,y,z;q) = q \sumstar_{\substack{a \shortmod{q} \\ a \neq - \overline{z}}} 
e\left(\frac{a x + a (\overline{1 + z a}) y}{q} \right)
- \sumstar_{a} \sumstar_{b} e\left(\frac{a x + b y }{q} \right),
\end{equation*}
where in case $q | z$ we interpret the condition $a \neq -\overline{z}$ to be vacuous.  If $q|z$ we get
\begin{equation*}
T(x,y,0;q) = qc_q(x+y) - c_q(x) c_q(y),
\end{equation*}
where $c_q(n)$ is the Ramanujan sum.
To continue the calculation we impose the condition $(z,q) = 1$.
Using $a \rightarrow \overline{z} a - \overline{z}$ and calculating the sum over $a$ as a Kloosterman sum, we get
\begin{equation*}
T(x,y,z;q) = q S(x\overline{z}, -y \overline{z};q) e\left(\frac{-x \overline{z} + y \overline{z}}{q}\right) - q - c_q(x) c_q(y).
\end{equation*}
Finally, we obtain \eqref{eq:Tbound} from case analysis, the Weil bound, and the easy computation $c_q(n) = q-1$ if $q|n$, and $-1$ if $q \nmid n$.
\end{proof}

\begin{proof}[Proof of Proposition \ref{thm:expsum}]
Applying Weil's bound for Kloosterman sums to the sum over $k$ (by completing the sum) gives
\begin{equation}
S(K,L;q) \ll L q^{\half + \varepsilon} + KLq^{-1} \ll L q^{\half + \varepsilon}.
\end{equation}
Further gains can be obtained by using the sum over $l$.
To begin, we smooth the sum $S$ to
\begin{equation*}
S_0(K, L) := \sumstar_{l}  \left| \sumstar_{k}  e\left(\frac{l \overline{k}}{q}\right) w\left(\frac{k}{K}\right) \right| w_0 \left(\frac{l}{L} \right),
\end{equation*}
where $w_0$ is a nonnegative, Schwartz-class function satisfying $w_0(x) \geq 1$ for $0\leq x \leq 1$.  By positivity, $S \leq S_0$.  Applying Cauchy's inequality gives
\begin{equation}
\label{eq:sweatshirt}
S_0(K, L)^2 \ll L    \sumstar_{k_1} \sumstar_{k_2} \sumstar_{l} e\left(\frac{l (\overline{k_1} - \overline{k_2}}{q}\right) w\left(\frac{k_1}{K}\right) w\left(\frac{k_2}{K}\right) w_0 \left(\frac{l}{L} \right).
\end{equation}

Now apply Poisson summation in each of $l$, $k_1$, and $k_2 \pmod{q}$ to give that the right hand side above is
\begin{equation*}
\frac{L^2 K^2}{q^3} \sum_x \sum_y \sum_z T(x,y,z;q) \widehat{w}\left(\frac{x}{q/K}\right) \widehat{w}\left(\frac{y}{q/K}\right) \widehat{w_0}\left(\frac{z}{q/L}\right).
\end{equation*}
From \eqref{eq:Tbound} and a case analysis, we get
\begin{equation*}
S_0(K, L)^2 \ll \frac{L^2 K^2}{q^3} \left(q^{\frac{3}{2}} \frac{q^3}{L K^2} + q^2 \frac{q}{K} \right)q^{\varepsilon}
= L q^{\frac{3}{2}+\varepsilon} + L^2 K q^{\varepsilon}.
\qedhere
\end{equation*}
\end{proof}

\subsection{An estimate for $N_1$ small}
\label{section:N1small}
In this section we prove Lemma \ref{lemma:chewyquaker}.
We continue with the estimation of \eqref{eq:brain}, recalling that the error terms prior to arriving at \eqref{eq:brain} are acceptable for Lemma \ref{lemma:chewyquaker}.  It is instructive to consider the crucial case $M = q^{\half}$, $N = q^{\frac32}$, $N_1 = q^{\half}$, $N_2 = q$ (whence $H = 1$).  We begin with the simple reciprocity law
\begin{equation*}
e\left(\frac{hm \overline{n_1}}{q}\right) = e\left(\frac{-hm \overline{q}}{n_1}\right)e\left(\frac{hm }{q n_1}\right),
\end{equation*}
where $q \overline{q} \equiv 1 \pmod{n_1}$.  Using $e(\frac{hm}{qn_1}) = 1 + O(\frac{hm}{qn_1})$ and $\frac{\phi(q)}{q \phi^*(q)} = q^{-1} + O(q^{-2})$ gives
\begin{equation*}
R = \frac{N_2}{q \sqrt{MN}} \sum_{h \neq 0} \widehat{W}\left(\frac{h}{H}\right) \sum_{n_1}  W\left(\frac{n_1}{N_1}\right)  \sum_{m} d(m) e\left(\frac{-hm \overline{q}}{n_1}\right) W\left(\frac{m}{M}\right)
+ O(\frac{M^{\frac32}}{N^{\frac32}} N_1 q^{\varepsilon}).
\end{equation*}
Let $R_1$ be the new sum above.  Let $g = (h, n_1)$, whence
\begin{equation*}
R_1 = \frac{N_2}{q \sqrt{MN}} \sum_{g \geq 1} \sum_{h \neq 0} \widehat{W}\left(\frac{h}{H/g}\right) \sum_{(n_1, h) = 1}  W\left(\frac{n_1}{N_1/g}\right) U(h,m, n_1),
\end{equation*}
\begin{equation*}
U(h,m, n_1) = \sum_{m} d(m) e\left(\frac{-hm \overline{q}}{n_1}\right) W\left(\frac{m}{M}\right).
\end{equation*}

Using the Mellin transform of $W$ to write the sum over $m$ in terms of the Estermann function gives
\begin{equation}
\label{eq:pizzaboat}
U(h,m, n_1)
=
\frac{1}{2 \pi i} \int_{(2)} M^{s}  D(s, 0, \frac{-h\overline{q}}{n_1}) \widetilde{W}(s) ds.
\end{equation}
Next move $s$ to $-1$, crossing a double pole at $s=1$ which gives a `main term' of size $\ll M n_1^{-1} q^{\varepsilon}$, which contributes $\ll M^{\half} N^{-\half} q^{\varepsilon}$ to $R_1$, an acceptable error term.  Applying the functional equation of the Estermann function (Lemma \ref{lem:Dq}) gives $U(h,m,n_1) = U_{+}(h,m,n_1) + U_{-}(h,m,n_1) + O(M n_1^{-1} q^{\varepsilon})$,
where
\begin{equation}
\label{eq:blackshoes}
U_{\pm}(h,m,n_1)= \frac{1}{2 \pi i} \int_{(-1)} M^{s} n_1^{1-2s}  D(1-s, 0, \frac{\mp q\overline{h}}{n_1}) \widetilde{W}_{\pm}(s) ds,
\end{equation}
where
\begin{equation*}
\widetilde{W}_{\pm}(s) = \widetilde{W}(s) 2 (2\pi)^{2s-2} \Gamma(1-s)^2 C_{\pm}(\pi s),
\end{equation*}
and $C_{+} = 1$, $C_{-} = - \cos$.  Expanding $D$ into absolutely convergent Dirichlet series and letting $s \rightarrow -s + 1$ gives that
\begin{equation*}
U_{\pm}(h,m,n_1) = \frac{M}{n_1} \sum_{k} d(k) e\left(\frac{\mp k q \overline{h}}{n_1} \right)  
y_{\pm}\left(\frac{kM}{n_1^2}\right),
\end{equation*}
where
\begin{equation*}
y_{\pm}(x) = \pm \frac{1 }{2 \pi i} \int_{(\varepsilon)} x^{-s} \widetilde{W}(1-s) 2 (2\pi)^{-s} \Gamma(s)^2 C_{\pm}(\pi s) ds.
\end{equation*}
Note $x^j y_{\pm}^{(j)}(x) \ll_A x^{-A}$ for any $A > 0$ (by taking $\text{Re}(s)$ large) and $y_{\pm}(x) = c_1 \log{x} + c_2 + O(x^{1-\varepsilon})$ for $x$ small (by taking $\text{Re}(s) = -1 + \varepsilon$).  Also $y_{\pm}^{(j)}(x) = c_j + O(x^{1-\varepsilon})$ for $j=1,2,\ldots$ and $x$ small.  

Thus $R_1 = R_{+} + R_{-} + O(M^{\half} N^{-\half} q^{\varepsilon})$, where
\begin{equation}
R_{\pm} = \frac{MN_2}{q N_1 \sqrt{MN}} \sum_g g \sum_{h \neq 0} \widehat{W}\left(\frac{hg}{H}\right) 
\sum_{(n_1, h) = 1} W_0\left(\frac{n_1g}{N_1}\right)
 \sum_{k} d(k) e\left(\frac{\mp k q \overline{h}}{n_1} \right)  
y_{\pm}\left(\frac{kM}{n_1^2}\right),
\end{equation}
where here $W_0(x) = x^{-1} W(x)$.  The trivial bound gives
\begin{equation}
R_{\pm} \ll \frac{N_1^2}{\sqrt{MN}} q^{\varepsilon}.
\end{equation}

It now suffices to show
\begin{equation}
\label{eq:Rpmbound}
R_{\pm} \ll \frac{N_1^2 q^{1/2 + \varepsilon}}{N^{5/6} M^{2/3}},
\end{equation}
which we proceed to prove.
Let $V(h,k)$ be given by
\begin{equation*}
V(h,k) = \sum_{(n_1, h) = 1} e\left(\frac{k q \overline{h}}{n_1} \right)  
y\left(\frac{kM}{n_1^2}\right) W_0\left(\frac{n_1}{N_1/g}\right).
\end{equation*}
Our next goal is to use the theory of exponential sums to prove that there is cancellation in the sum over $n_1$, and bound the sums over $h$ and $k$ trivially.  The presence of $\overline{h}$ is a bit of a nuisance, so we again use the elementary reciprocity to law to write
\begin{equation*}
e\left(\frac{k q \overline{h}}{n_1} \right) = e\left(\frac{-k q \overline{n_1}}{h} \right) e\left(\frac{k q}{h n_1} \right).
\end{equation*}
By splitting $n_1$ into residue classes $\pmod{h}$, we get
\begin{equation*}
V(h,k) = \sumstar_{a \shortmod{h}} e\left(\frac{-k q \overline{a}}{h} \right) \sum_{r}  e\left(\frac{k q}{h (a + hr)} \right)  
y\left(\frac{kM}{(a + hr)^2}\right) W\left(\frac{a + hr}{N_1/g}\right).
\end{equation*}
We may assume $gh \ll H q^{\varepsilon} = q N_2^{-1} q^{\varepsilon}$, which is $\ll N_1 q^{-\varepsilon}$ since our results are trivial if $N \ll q^{1 + \varepsilon}$ (the point is that there is room for summation over $r$).

By partial summation it suffices to bound
\begin{equation*}
\sum_{\frac{N_1}{gh} < r \leq \frac{N_1}{gh} + C} e\left(\frac{k q}{h (a + hr)} \right)
\end{equation*}
for $C \ll \frac{N_1}{gh}$.  We apply the convenient bound given by Corollary 8.5 of \cite{IK}, which we reproduce here (with slight changes in notation).
\begin{myprop}
Suppose $f(x)$ satisfies
\begin{equation*}
\frac{F}{A} \leq \frac{x^k}{k!} |f^{(k)}(x)| \leq F
\end{equation*}
for $k=2,3$ in the segment $[Q,2Q]$.  Then for $1 \leq Q' \leq Q \leq F$ we have
\begin{equation*}
\sum_{Q < m \leq Q'} e(f(m)) \ll A F^{\frac16} Q^{\half} \log{3Q},
\end{equation*}
where the implied constant is absolute.
\end{myprop}

For our application, 
\begin{equation*}
f(x) = \frac{kq}{h(a + hx)}, \quad F \asymp \frac{kqg}{N_1 h}, \quad Q = \frac{N_1}{gh},
\end{equation*} 
and $A$ is absolute.
Then we get
\begin{equation*}
V(h,k) \ll h \left(\frac{kgq}{N_1 h}\right)^{\frac16} \left(\frac{N_1}{gh}\right)^{\half} q^{\varepsilon}.
\end{equation*}
Using this bound gives
\begin{align}
R_{\pm} \ll \frac{MN_2}{q N_1 \sqrt{MN}} \left(\frac{q}{N_2}\right)^2 \frac{N_1^2}{M} \left(\frac{N}{M}\right)^{\frac16} \left(\frac{N}{q}\right)^{\half} q^{\varepsilon},
\end{align}
which simplifies to give \eqref{eq:Rpmbound}, as desired.  This completes the proof of Lemma \ref{lemma:chewyquaker}.

\section{A binary divisor sum}
\label{section:binary}
When $M$ and $N$ are roughly the same size we treat $B_{M,N}$ as a binary divisor sum.  We loosely follow the presentation of Motohashi, Chapter 4 \cite{M}.  Before beginning the treatment of the divisor sum, we first present some material to which we shall shortly refer.

\subsection{Some arithmetical sums}
We require the computation of some arithmetical sums.
\begin{mylemma}
\label{lem:divisorfunctionsum}
Let $\text{Re}(s - \lambda) > 1$ and $\text{Re}(s) > 1$.  Then if $d$ is either $1$ or a prime then
\begin{equation}
\sum_{n \equiv 0 \shortmod{d}} \frac{\sigma_{\lambda}(n)}{n^s} =  \zeta(s) \zeta(s - \lambda) d^{-s} \sum_{\substack{bc = d}} \frac{\mu(b)}{b^{s - \lambda}} \sigma_{\lambda}(c).
\end{equation}
\end{mylemma}
Actually Lemma \ref{lem:divisorfunctionsum} holds for general $d$ but we only need it for $d =1$ or prime.  In this special case the proof is greatly simplified so we omit it.

\begin{mylemma}
If $\text{Re}(\alpha) < 0$ we have
\begin{equation}
\label{eq:ramanujandivisor}
\sigma_{\alpha}(n) = \zeta(1 - \alpha) \sum_{l} \frac{c_l(n)}{l^{1 - \alpha}}.
\end{equation}
\end{mylemma}
\begin{proof}
This follows immediately from using the sum of divisors formula for the Ramanujan sum, that is $c_l(n) = \sum_{d|(l,n)} d \mu(l/d)$, and reversing orders of summation.
\end{proof}
\begin{mylemma}
\label{prop:handlsum}
Let $d$ be either $1$ or prime, $\text{Re}(s) > 1$, and $\text{Re}(\lambda) > -1$. Then
\begin{equation}
\sum_{l} \frac{1}{l^{2+\lambda}}
\sumstar_{h \shortmod{l}}  \sum_{n \equiv 0 \shortmod{d}} \frac{e(\frac{nh}{l})}{n^s}
=
\frac{\zeta(s)\zeta(1 + \lambda + s)}{d^s \zeta(2+ \lambda)} \left(1 + \frac{1}{d^{1 + \lambda}} - \frac{1}{d^{1+\lambda + s}} \right).
\end{equation}
\end{mylemma}
\begin{proof}
Evaluating the sum over $h$ gives that the left hand side above equals
\begin{equation}
\sum_{n \equiv 0 \shortmod{d}} \sum_{l} \frac{1}{l^{2+ \lambda}} \frac{ c_l(n)}{n^s}.
\end{equation}
Executing the sum over $l$ using \eqref{eq:ramanujandivisor} shows it equals
\begin{equation}
\frac{1}{\zeta(2+ \lambda)} \sum_{n \equiv 0 \shortmod{d}} \frac{\sigma_{-1 - \lambda}(n)}{ n^s}.
\end{equation}
Finally using Lemma \ref{lem:divisorfunctionsum} and the fact that $d$ is either $1$ or a prime completes the proof.
\end{proof}
\subsection{An approximate functional equation for the divisor function}
\begin{mylemma}
\label{lemma:ramanujan}
For any positive integer $n$ and $\lambda \in \mc$, 
\begin{equation}
\sigma_{\lambda}(n) = \sum_{l} \frac{c_l(n)}{l^{1 - \lambda}}  f_{\lambda}(\frac{l}{\sqrt{n}})
+ n^{\lambda} \sum_{l} \frac{c_l(n)}{l^{1 + \lambda}}  f_{-\lambda}(\frac{l}{\sqrt{n}}),
\end{equation}
where
\begin{equation}
f_{\lambda}(x) = \int_{(a)} x^{-w} \zeta(1 - \lambda + w) \frac{G(w)}{w} dw,
\end{equation}
$a > |\text{Re}(\lambda)|$, $c_l(n)$ is the Ramanujan sum, and $G$ is as in Definition \ref{def:G}.
\end{mylemma}
Remarks.
\begin{itemize}
\item The reason to expand the divisor function into such a series is that the exponential sum formula for the Ramanujan sum $c_l(n+f)$ will allow for the separation of variables $n$ and $f$ in $\sigma_{\lambda}(n+f)$.  This is a simple alternative to the delta method \cite{DFI}.
\item Motohashi used the formula \eqref{eq:ramanujandivisor} to accomplish the separation of variables.  However, \eqref{eq:ramanujandivisor} does not hold for $\alpha$ in a neighborhood of the origin.
\item An approximate functional equation is essentially equivalent to a functional equation.  The sum of divisors function $\sigma_{\lambda}(n)$ satisfies the functional equation
\begin{equation*}
\sigma_{\lambda}(n) = n^{\lambda} \sigma_{-\lambda}(n).
\end{equation*}
The function $n^{-\frac{\lambda}{2}} \sigma_{\lambda}(n)$ is perhaps more natural to study because it is invariant under $\lambda \leftrightarrow -\lambda$.  Of course, these appear as Fourier coefficients of Eisenstein series \eqref{eq:fouriereisenstein}.
\end{itemize}
\begin{proof}
By a contour shift we have
\begin{equation*}
 \int_{(a)} \sigma_{\lambda - w}(n) n^{w/2}\frac{G(w)}{w} dw = \sigma_{\lambda}(n) + \int_{(-a)} \sigma_{\lambda - w}(n) n^{w/2} \frac{G(w)}{w} dw.
\end{equation*}
An application of the functional equation $\sigma_{\lambda - w}(n) = n^{\lambda -w} \sigma_{-\lambda + w}(n)$ and the change of variables $w \rightarrow -w$ gives
\begin{equation*}
\sigma_{\lambda}(n) = \int_{(a)} \sigma_{\lambda - w}(n) n^{w/2} \frac{G(w)}{w} dw
+ n^{\lambda}\int_{(a)} \sigma_{-\lambda - w}(n) n^{w/2} \frac{G(w)}{w} dw.
\end{equation*}
Inserting \eqref{eq:ramanujandivisor} into the above integral representation and reversing the order of summation and integration completes the proof.
\end{proof}
\subsection{Separation of variables}
We begin the treatment of $B_{M,N}$ by solving the congruence $m \equiv n \pmod{d}$ by letting $n = m + f$, where $d|f$, and $f > 0$.  We have
\begin{multline*}
B_{M, N} =  \frac{1}{\phi^{*}(q)} \sum_{d | q} \phi(d) \mu(\frac{q}{d}) 
\sum_{\substack{f \equiv 0 \shortmod{d}}}{}
\\
\sum_{m}  \frac{\sigma_{\alpha - \beta}(m) \sigma_{\gamma - \delta}(m+f)}{m^{\half + \alpha} (m+f)^{\half + \gamma} } V_{\alpha, \beta, \gamma, \delta}\left(\frac{m(m+f)}{q^2}\right) W\left(\frac{m}{M}, \frac{m+f}{N}\right).
\end{multline*}
Our immediate goal is to separate the variables $n$ and $f$ both arithmetically and analytically.
We use the expansion of $\sigma_{\gamma - \delta}(n+f)$ into Ramanujan series given by Lemma \ref{lemma:ramanujan} to arithmetically separate the variables $n$ and $f$.  Let $C_{M, N}$ be the contribution from the first term of Lemma \ref{lemma:ramanujan} and $\widetilde{C}_{M, N}$ be the second term.  Then
\begin{equation}
\label{eq:BC}
B_{M,N} = C_{M,N} + \widetilde{C}_{M, N}, \quad \widetilde{C}_{M,N}(\alpha, \beta, \gamma, \delta) = C_{M,N}(\alpha, \beta, \delta, \gamma),
\end{equation}
where the formula for $C_{M,N}$ is
\begin{multline}
\label{eq:cupcake}
C_{M,N} = \frac{1}{\phi^{*}(q)} \sum_{d | q} \phi(d) \mu(\frac{q}{d}) \sum_{l} \frac{1}{l^{1 - \gamma + \delta}} 
\sum_{\substack{f \equiv 0 \shortmod{d} }}{}
\\
\sum_{m}  \frac{\sigma_{\alpha - \beta}(m) c_l(m+ f)}{m^{\half + \alpha} (m+f)^{\half + \gamma} } V_{\alpha, \beta, \gamma, \delta}\left(\frac{m(m+f)}{q^2}\right) W\left(\frac{m}{M}, \frac{m+f}{N}\right) f_{\gamma-  \delta}\left(\frac{l}{\sqrt{m+f}}\right)
\end{multline}
The function $f_{\gamma-\delta}$ should not be confused with the variable of summation $f$; at any rate, this alphabetical accident shall clear up shortly.  Before stating our preferred formula for $C_{M,N}$ we set some notation.
Define the Mellin pair
\begin{equation*}
\begin{cases}
\widetilde{W}(u_1, u_2) = \int_0^{\infty} \int_0^{\infty} W(x,y) x^{u_1} y^{u_2} \frac{dx dy }{xy}, \\
W(x,y)  = \left(\frac{1}{2 \pi i}\right)^2 \int_{(c_{u_2})} \int_{(c_{u_1})} \widetilde{W}(u_1, u_2) x^{-u_1} y^{-u_2} du_1 du_2.
\end{cases}
\end{equation*}
We generally use the notation $\int_{(c_s)}$ to denote the vertical line of integration with $\text{Re}(s) = c_s$, and similarly for other subscripts.  Let
\begin{equation}
\label{eq:Hthree}
H_1(s, u_1, u_2, w) = \frac{G(s) G(w)}{s w} g(s) \widetilde{W}(u_1, u_2) \zeta(1 - \gamma + \delta + w),
\end{equation}
\begin{equation}
H(s, u_1, u_2, v, w) = \frac{\Gamma(v) \Gamma(\half + \gamma + s  + u_2 -v -\frac{w}{2})}{\Gamma(\half + \gamma + s  + u_2  -\frac{w}{2})} H_1(s,u_1,u_2,w).
\end{equation}
We claim both $H_1$ and $H$ have rapid decay as any of the variables get large in the imaginary direction.  This is easy to see for $H_1$. For $H$ we note that a crude use of Stirling's approximation shows that (for $b$ and $v$ in some fixed vertical strip) 
\begin{equation}
\label{eq:gammavb}
\frac{\Gamma(v) \Gamma(b-v)}{\Gamma(b)} \ll 
(1+|v|)^{-A} (1+|b|)^{2A}, 
\end{equation}
where $A$ is any sufficiently large (depending on the fixed vertical strip) positive number, and for $b$ and $v$ avoiding the poles of the gamma function.  For example, if $|v| \leq 2|b|$ then the ratio of gammas has exponential decay in $v$, uniformly in $b$.  In the opposite case just use the fact that it is bounded by a polynomial in $v$ and $b$ and multiply and divide by $(1+|v|)^{-A}$.  The rapid decay of $H_1$ in all variables overcomes this potential growth in $b$ to show $H$ has rapid decay in all variables.
\begin{mylemma}
\label{lemma:POintegral}
With $c_s=2$, $c_v = c_w = \varepsilon$, and $c_{u_1} = c_{u_2} = 0$, we have
\begin{multline}
\label{eq:POintegral}
C_{M,N} = \frac{1}{\phi^{*}(q)}  \sum_{d | q} \phi(d) \mu(\frac{q}{d}) \sum_{l} \frac{1}{l^{1 - \gamma + \delta}}
\sumstar_{h \shortmod{l}}    
\sum_{\substack{f \equiv 0 \shortmod{d}}}{} \frac{e\left(\frac{hf}{l} \right) }{f^{\half + \gamma}}
\left(\frac{1}{2 \pi i} \right)^5 \int_{(c_s)} \int_{(c_w)}
\\
 \int_{(c_{u_2})} \int_{(c_{u_1})} \int_{(c_v)}   \frac{q^{2s} M^{u_1} N^{u_2} }{f^{s+ u_2 - v - \frac{w}{2} } l^w  }
D(\tfrac12 + \alpha + s + u_1 +v, \alpha - \beta, \frac{h}{l})
H(s, u_1, u_2, v, w)
dv du_1 du_2 dw ds.
\end{multline}
\end{mylemma}

\begin{proof}
Inserting the exponential sum formula for $c_l(m+f)$ into \eqref{eq:cupcake} gives
\begin{multline*}
C_{M,N} = \frac{1}{\phi^{*}(q)} \sum_{d | q} \phi(d) \mu(\frac{q}{d}) \sum_{l} \frac{1}{l^{1 - \gamma + \delta}}
\sumstar_{h \shortmod{l}} 
\sum_{\substack{f \equiv 0 \shortmod{d}}}{} e\left(\frac{hf}{l}\right) 
\\
\sum_{m}  \frac{\sigma_{\alpha - \beta}(m) e\left(\frac{mh}{l}\right)}{m^{\half + \alpha} (m+f)^{\half + \gamma} } V_{\alpha, \beta, \gamma, \delta}\left(\frac{m(m+f)}{q^2}\right) W\left(\frac{m}{M}, \frac{m+f}{N}\right) f_{\gamma - \delta}\left(\frac{l}{\sqrt{m+f}}\right).
\end{multline*}
Next we separate the variables $m$ and $f$ analytically by taking Mellin transforms.  We have
\begin{multline}
\label{eq:analyticsep}
\frac{1}{(m + f)^{\half + \gamma}} V\left(\frac{m(m+f)}{q^2}\right) f_{\gamma - \delta}\left(\frac{l}{\sqrt{m+f}}\right) W\left(\frac{m}{M}, \frac{m+f}{N}\right) 
\\ = \left(\frac{1}{2 \pi i} \right)^4 \int_{(c_s)} \int_{(c_w)} \int_{(c_{u_1})} \int_{(c_{u_2})}   \frac{q^{2s} M^{u_1} N^{u_2} H_1(s, u_1, u_2, w)}{l^w m^{s + u_1}  (m+f)^{\half + \gamma +s +u_2 - \frac{w}{2}}}du_2 du_1 dw ds.
\end{multline}
To separate $m$ and $f$ in $(m + f)^{-\half - \gamma - s -u_2 + \frac{w}{2}}$ we use the following formula (17.43.21 of \cite{GR}):
\begin{align}
(1 + x)^{-b}
=
\frac{1}{2 \pi i} \int_{(c_v)} \frac{\Gamma(v) \Gamma(b-v)}{ \Gamma(b)} x^{-v} dv,
\end{align}
valid for $0 < c_v < \text{Re}(b)$.   Recall \eqref{eq:gammavb} for the convergence.

Thus \eqref{eq:analyticsep} equals
\begin{equation}
\label{eq:Canalyticsep}
\left(\frac{1}{2 \pi i} \right)^5 \int_{(c_s)} \int_{(c_w)} \int_{(c_{u_2})} \int_{(c_{u_1})} \int_{(c_v)}  \frac{q^{2s}  M^{u_1} N^{u_2} H(s, u_1, u_2, v, w)}{f^{\half + \gamma + s+ u_2 - v - \frac{w}{2}} l^w m^{s + u_1 + v} } dv du_1 du_2 dw ds.
\end{equation}
Inserting \eqref{eq:Canalyticsep} into $C_{M,N}$, we obtain
\begin{multline}
\label{eq:blackvitamin}
C_{M,N} = \frac{1}{\phi^{*}(q)}  \sum_{d | q} \phi(d) \mu(\frac{q}{d}) \sum_{l} \frac{1}{l^{1 - \gamma + \delta}}
\sumstar_{h \shortmod{l}}    
\sum_{\substack{f \equiv 0 \shortmod{d}}}{} \frac{e\left(\frac{hf}{l} \right) }{f^{\half + \gamma}}
\sum_{m}  \frac{\sigma_{\alpha - \beta}(m) e\left(\frac{mh}{l}\right)}{m^{\half + \alpha}} 
\\
\left(\frac{1}{2 \pi i} \right)^5 \int_{(c_s)} \int_{(c_w)} \int_{(c_{u_2})} \int_{(c_{u_1})} \int_{(c_v)}   \frac{q^{2s} M^{u_1} N^{u_2} H(s, u_1, u_2, v, w)}{f^{ s+ u_2 - v - \frac{w}{2}} l^w m^{s + u_1 + v}  }dv du_1 du_2 dw ds.
\end{multline}
The sum over $m$ converges absolutely on the contours of integration as stated in Lemma \ref{lemma:POintegral}.  Writing this sum over $m$ in terms of the Estermann $D$-function finishes the proof.
\end{proof}

Our goal is to develop the asymptotics of $C_{M,N}$ by moving the lines of integration past the poles of the $D$-function (which will contribute the main terms) and to analyze the remainder term using estimates for sums of Kloosterman sums using the Kuznetsov formula.

To begin, we move the $s$-line of integration to $c_{s} = \varepsilon$, passing the two poles of $D$ at $\half + s + u_1 + v = 1- \alpha$ and $\half + s + u_1 + v = 1- \beta$.  Let $E_{M,N}(\alpha,\beta,\gamma,\delta)$ be the integral on the new line of integration, and let $P_{M,N} = P_{M,N}(\alpha,\beta,\gamma,\delta)$ be the contribution of the former pole; by symmetry considerations, the latter pole is $P_{M,N}(\beta,\alpha,\gamma,\delta)$.  Thus we have
\begin{equation}
\label{eq:CPE}
C_{M,N}(\alpha,\beta,\gamma,\delta) = P_{M,N}(\alpha,\beta,\gamma,\delta) + P_{M,N}(\beta,\alpha,\gamma,\delta) + E_{M,N}(\alpha,\beta,\gamma,\delta).
\end{equation}
We bound $E_{M,N}$ in Section \ref{section:errorterms}, proving Theorem \ref{thm:Ebound}, and continue with our calculation of $P_{M,N}$ in Section \ref{section:mainterms}.

\section{The main terms}
\label{section:mainterms}
The partition of unity is an obstruction in the computation of these main terms.  It turns out to be easier to sum over $M$ and $N$ before doing finer analysis of the main terms.  
\begin{mylemma}
\label{lemma:claimedupperboundonP}
For $N \gg M$, we have
\begin{equation}
\label{eq:claimedupperboundonP}
P_{M,N} \ll N^{-\half} M^{\half} q^{\varepsilon}.
\end{equation}
\end{mylemma}
We prove Lemma \ref{lemma:claimedupperboundonP} in Section \ref{section:upperboundonP}.
\subsection{Recomposition}
Define $P := \sum_{M,N} P_{M,N}$.
\begin{mylemma}
\label{lemma:POmainterm}
We have
\begin{multline}
\label{eq:POmainterm}
P(\alpha, \beta, \gamma, \delta) =
\zeta(1 - \alpha + \beta)    
\left(\frac{1}{2 \pi i} \right)^2 
\int_{(\varepsilon)} \int_{(\half -2\varepsilon)} 
 q^{-\alpha - \gamma + \frac{w}{2}}
\frac{G(v)G(w) g(v)}{vw}
\\
\frac{\Gamma(\half - \alpha -v) \Gamma(\alpha + \gamma +2v - \frac{w}{2})}{\Gamma(\half + \gamma + v - \frac{w}{2})}
\zeta(1 - \gamma + \delta + w) 
\frac{\zeta( \alpha + \gamma +2v - \frac{w}{2}) \zeta(1 + \beta + \delta +2v + \frac{w}{2})}{\zeta(2-\alpha + \beta -\gamma + \delta + w)}
dv dw 
\\
+ O(q^{-1+\varepsilon}).
\end{multline}
\end{mylemma} 
The four main terms from $B_O'$ are obtained by taking the main terms of $B_O$ and switching $\alpha$ and $\gamma$, and $\beta$ and $\delta$.  One term is
\begin{multline}
\label{eq:mt2}
\zeta(1 - \gamma + \delta)    
\left(\frac{1}{2 \pi i} \right)^2 
\int_{(\varepsilon)} \int_{(\half -2\varepsilon)} 
 q^{-\alpha - \gamma + \frac{w}{2}} 
\frac{G(v)G(w) g(v)}{vw}
\zeta(1 - \alpha + \beta + w) 
\\
\frac{\Gamma(\half - \gamma -v) \Gamma(\alpha + \gamma +2v - \frac{w}{2})}{\Gamma(\half + \alpha + v - \frac{w}{2})}
\frac{\zeta( \alpha + \gamma +2v - \frac{w}{2}) \zeta(1 + \beta + \delta +2v + \frac{w}{2})}{\zeta(2-\alpha + \beta -\gamma + \delta + w)}
dv dw,
\end{multline}
and similar formulas hold for the other terms.
\begin{proof}
We first show how to evaluate the sum over $M$ and $N$.  As a warmup problem, we show that for a `nice' function $F$, 
\begin{equation*}
\sum_{M, N} \left(\frac{1}{2 \pi i} \right)^2  \int_{(c_{u_2})} \int_{(c_{u_1})} F(u_1, u_2) \widetilde{W}_{M,N}(u_1, u_2) du_1 du_2 = F(0,0),
\end{equation*}
For a proof, let $f$ be the inverse Mellin transform of $F$ and begin with the Mellin convolution formula to write the left hand side above as
\begin{equation*}
\sum_{M, N} \int_0^{\infty} \int_0^{\infty} f(x,y) W_{M,N} (x^{-1}, y^{-1}) \frac{dx}{x} \frac{dy}{y} = 
\int_0^{\infty} \int_0^{\infty} f(x,y) \frac{dx}{x} \frac{dy}{y} =F(0,0).
\end{equation*}
The point is that on the level of the Mellin transform, we divide out by $\widetilde{W}(0,0)$ and evaluate everything at $u_1 = u_2 =0$.

The computation of $P$ is of this form but with the extra step of computing the
residue of the Estermann function at $\half + s + v = 1 - \alpha$.  We make a brief diversion to justify this step.

Note that $P_{M,N}$ is the difference of two integrals of the form \eqref{eq:POintegral} with the only difference between the two being the placement of the line of integration $c_s$.  Applying the changes of variable $s \rightarrow s - u_1$ and $w \rightarrow w - 2u_1 + 2u_2$ transforms \eqref{eq:POintegral} into
\begin{multline}
\label{eq:hat}
\frac{1}{\phi^{*}(q)}  \sum_{d | q} \phi(d) \mu(\frac{q}{d}) \sum_{l} \frac{1}{l^{1 - \gamma + \delta}}
\sumstar_{h \shortmod{l}}    
\sum_{\substack{f \equiv 0 \shortmod{d}}}{} \frac{e\left(\frac{hf}{l} \right) }{f^{\half + \gamma}}
\left(\frac{1}{2 \pi i} \right)^5 \int_{(c_s)} \int_{(c_w)}
 \int_{(c_{u_2})} 
\\
\int_{(c_{u_1})} \int_{(c_v)}   \frac{q^{2s -2u_1} M^{u_1} N^{u_2} }{f^{s - v - \frac{w}{2} } l^{w - 2u_1 + 2u_2}  }
D(\tfrac12 + \alpha + s  +v, \alpha - \beta, \frac{h}{l})
H_2(s, u_1, u_2, v, w)
dv du_1 du_2 dw ds,clean 
\end{multline}
where
\begin{multline*}
H_2(s, u_1, u_2, v, w) = \frac{\Gamma(v) \Gamma(\half + \gamma + s -v- \frac{w}{2})}{\Gamma(\half + \gamma + s - \frac{w}{2})} 
\frac{G(s - u_1) G(w - 2u_1 + 2u_2)}{(s- u_1) (w - 2u_1 + 2u_2)} 
\\
g(s-u_1) \widetilde{W}(u_1, u_2) \zeta(1 - \gamma + \delta + w - 2u_1 + 2u_2).
\end{multline*}
Hence after summing over $M$ and $N$ in \eqref{eq:hat} we get
\begin{multline}
\label{eq:hat2}
\frac{1}{\phi^{*}(q)}  \sum_{d | q} \phi(d) \mu(\frac{q}{d}) \sum_{l} \frac{1}{l^{1 - \gamma + \delta}}
\sumstar_{h \shortmod{l}}    
\sum_{\substack{f \equiv 0 \shortmod{d}}}{} \frac{e\left(\frac{hf}{l} \right) }{f^{\half + \gamma}}
\\
\left(\frac{1}{2 \pi i} \right)^3 \int_{(c_s)} \int_{(c_w)} \int_{(c_v)}   \frac{q^{2s} }{f^{s - v - \frac{w}{2} } l^w  }
D(\tfrac12 + \alpha + s  +v, \alpha - \beta, \frac{h}{l})
H_3(s, v, w)
dv  dw ds,
\end{multline}
where $H_3$ is the following function with rapid decay in all variables (recall \eqref{eq:gammavb}):
\begin{equation}
H_3(s,v,w) = \frac{\Gamma(v) \Gamma(\half + \gamma + s -v- \frac{w}{2})}{\Gamma(\half + \gamma + s - \frac{w}{2})} 
\frac{G(s) G(w)}{s w} g(s) \zeta(1 - \gamma + \delta + w).
\end{equation}
Thus $P(\alpha,\beta,\gamma,\delta)$ is the residue of the integrand in \eqref{eq:hat2} at $s + v = \half - \alpha$.

Using Lemma \ref{lem:Dq}, we find that the residue of the Estermann function at this point is 
$ l^{-1+\alpha-\beta} \zeta(1 - \alpha + \beta)$,
and hence
\begin{multline*}
P=\frac{\zeta(1 - \alpha + \beta)}{\phi^{*}(q)} \sum_{d | q} \phi(d) \mu(q/d) \sum_{l} \frac{1}{l^{2 - \alpha + \beta - \gamma + \delta}}
\sumstar_{h \shortmod{l}}   
\sum_{\substack{f \equiv 0 \shortmod{d}}}{} e\left(\frac{hf}{l} \right)
\\
\left(\frac{1}{2 \pi i} \right)^2  \int_{(c_w)} \int_{(c_v)}   \frac{q^{1-2\alpha - 2 v} }{f^{1 - \alpha + \gamma - 2v - \frac{w}{2} } l^w  }
H_3(\tfrac12 - \alpha - v, v, w)
dv dw.
\end{multline*}
Apply the change of variables $v \rightarrow \half - \alpha - v$ to get
\begin{multline}
\label{eq:spoon}
P=\frac{\zeta(1 - \alpha + \beta)}{\phi^{*}(q)} \sum_{d | q} \phi(d) \mu(q/d) \sum_{l} \frac{1}{l^{2 - \alpha + \beta - \gamma + \delta}}
\sumstar_{h \shortmod{l}}   
\sum_{\substack{f \equiv 0 \shortmod{d}}}{} e\left(\frac{hf}{l} \right)
\\
\left(\frac{1}{2 \pi i} \right)^2  \int_{(c_w)} \int_{(c_v)}   \frac{q^{2 v} }{f^{ \alpha + \gamma  + 2v - \frac{w}{2} } l^w  }
H_3(v, \tfrac12 - \alpha - v, w)
dv dw.
\end{multline}
Note that the ratio of gamma functions in $H_3$ above is
\begin{equation*}
\frac{\Gamma(\half - \alpha - v) \Gamma( \alpha + \gamma + 2v - \frac{w}{2})}{\Gamma(\half + \gamma + v - \frac{w}{2})}.
\end{equation*}
Since $G(\half - \alpha) = 0$ (recall Definition \ref{def:G}), we may move $c_v$ to $1$ without encountering any poles.  The summation over $f$ converges absolutely, and we may execute the summations over $f$, $h$, and $l$ using Lemma \ref{prop:handlsum} to obtain that
\begin{multline}
\label{eq:messyd}
P= \frac{\zeta(1 - \alpha + \beta)}{\phi^{*}(q)} \sum_{d | q} \phi(d) \mu(q/d) 
\left(\frac{1}{2 \pi i} \right)^2  \int_{(c_w)} \int_{(c_v)} q^{2v}  H_3(v, \tfrac12 - \alpha - v, w)
\\
\frac{\zeta(\alpha + \gamma +2v-\frac{w}{2}) \zeta(1 + \beta + \delta + 2v + \frac{w}{2})}{d^{\alpha + \gamma +2v - \frac{w}{2}} \zeta(2-\alpha+\beta-\gamma+\delta + w)}
\left(1 + d^{-1 + \alpha - \beta + \gamma - \delta - w} - d^{-1 - \beta - \delta -2v- \frac{w}{2}} \right)
 dv  dw.
\end{multline}
Now we move the $v$-line of integration back to $c_v = \varepsilon$.  The pole at $\alpha + \gamma + 2v - \frac{w}{2} = 1$ gives
\begin{multline}
\label{eq:battery}
\frac{\zeta(1 - \alpha + \beta)}{\phi^{*}(q)} \sum_{d | q} \frac{\phi(d)}{d} \mu(q/d) 
\frac{1}{2 \pi i}  \int_{(\varepsilon)}  q^{-\alpha - \gamma + 1 + \frac{w}{2}}  
\\
\left(1 + d^{-1 + \alpha - \beta + \gamma - \delta - w} - d^{-2 +\alpha - \beta + \gamma - \delta - w} \right)
H_3(\tfrac12(1- \alpha - \gamma + \frac{w}{2}), \tfrac12(-\alpha + \gamma - \frac{w}{2}), w)
dw,
\end{multline}
which is seen to be $O(q^{-1 + \varepsilon})$ due to cancellation in the arithmetical sum over $d$.
Using trivial estimations shows that
\begin{multline*}
P= \zeta(1 - \alpha + \beta)
\left(\frac{1}{2 \pi i} \right)^2  \int_{(c_w)} \int_{(c_v)} q^{-\alpha - \gamma + \frac{w}{2}} 
\\
\frac{\zeta(\alpha + \gamma +2v-\frac{w}{2}) \zeta(1 + \beta + \delta + 2v + \frac{w}{2})}{\zeta(2-\alpha+\beta-\gamma+\delta + w)}
H_3(v, \tfrac12 - \alpha - v, w) + O(q^{-1+\varepsilon}).
\end{multline*}
Observe that
\begin{multline*}
H_3(v, \tfrac12 - \alpha - v, w)
=
\frac{G(v) G(w)}{v w} g(v) \zeta(1 - \gamma + \delta+ w) 
\frac{\Gamma(\tfrac12 - \alpha -v) \Gamma(\alpha + \gamma +2v - \frac{w}{2})}{\Gamma(\tfrac12 + \gamma + v - \frac{w}{2})},
\end{multline*}
which completes the proof.
\end{proof}

\subsection{An upper bound for $P_{M,N}$}
\label{section:upperboundonP}
In this section we prove Lemma \ref{lemma:claimedupperboundonP}.  A minor variation of some calculations in the proof of Lemma \ref{lemma:POmainterm} reduces the problem of bounding $P_{M,N}$ to bounding
\begin{multline*}
\zeta(1 - \alpha + \beta)
\left(\frac{1}{2 \pi i} \right)^4  \int_{(c_w)} \int_{(c_v)} \int_{(c_{u_1})} \int_{(c_{u_2})} q^{-\alpha - \gamma -u_1 -u_2+ \frac{w}{2}}  M^{u_1} N^{u_2}  \widetilde{W}(u_1, u_2)
\\
\frac{\zeta(\alpha + \gamma +u_1 + u_2 + 2v-\frac{w}{2}) \zeta(1 + \beta + \delta + u_1 + u_2 + 2v + \frac{w}{2})}{\zeta(2-\alpha+\beta-\gamma+\delta + w)}
\frac{G(v) G(w)}{v w} g(v)
\\
\frac{\Gamma(\half - \alpha - u_1 - v) \Gamma(\alpha + \gamma +u_1 + u_2 + 2v-\frac{w}{2})}{\Gamma(\half + \gamma + v  + u_2  -\frac{w}{2})}
  \zeta(1 - \gamma + \delta + w)
 du_2 du_1 dv dw,
\end{multline*}
where $c_{u_1} = c_{u_2} = 0$, and $c_v = c_w = \varepsilon$.  Then move $c_{u_1}$ to $\half - 2 \varepsilon$ followed by $c_{u_2} \rightarrow -\half$, and bound the integrand trivially with absolute values to finish the proof.

\section{The dual off-diagonal terms}
\label{section:dual}
To properly manipulate the main terms obtained in Section \ref{section:mainterms}, it is necessary to combine them with the corresponding main terms of $A_{\overline{O}}$.  The computation is similar to that of $B_{O}$ but there are some differences.   

As before, we require different methods of estimation depending on how close $m$ and $n$ are.  We apply the same partition of unity as before and write $A_{\overline{O}} = \sum_{M, N} A_{M,N}$ accordingly.  Clearly $A_{M,N}$ satisfies the bound \eqref{eq:faraway}.  

\subsection{Separation of variables}
We follow the methods of Section \ref{section:binary}.  Write $f = m + n$ so
\begin{multline*}
A_{M,N} =  \frac{1}{\phi^{*}(q)} \sum_{d | q} \phi(d) \mu(q/d) 
\sum_{f \equiv 0 \shortmod{d}} 
\\
\sum_{\substack{0 < m < f}} \frac{\sigma_{\alpha - \beta}(m) \sigma_{\gamma - \delta}(f-m)}{m^{\half + \alpha} (f-m)^{\half + \gamma} } V\left(\frac{m(f-m)}{q^2}\right)
W\left(\frac{m}{M}, \frac{f-m}{N}\right).
\end{multline*}
We use the expansion of $\sigma_{\gamma - \delta}(f-m)$ into Ramanujan series (Lemma \ref{lemma:ramanujan}) to separate the variables.  Write $A_{M,N} = C_{\overline{M,N}} + \widetilde{C}_{\overline{M,N}}$ to correspond to the first and second parts of the approximate functional equation.  Then
\begin{multline*}
C_{\overline{M,N}}  = \frac{1}{\phi^{*}(q)} \sum_{d | q} \phi(d) \mu(\frac{q}{d}) \sum_{l} \frac{1}{l^{1 - \gamma + \delta}} \sumstar_{h \shortmod{l}}\sum_{\substack{f \equiv 0 \shortmod{d} \\ f > 0}} e\left(\frac{-hf}{l}\right)
\\
\sum_{0<m<f}   \frac{\sigma_{\alpha - \beta}(m) e\left(\frac{hm}{l}\right) }{m^{\half + \alpha} (f -m)^{\half + \gamma} } V \left(\frac{m(f-m)}{q^2}\right) W\left(\frac{m}{M}, \frac{f-m}{N}\right).
\end{multline*}
Using the same methods as in the computation of $C_{M,N}$, we obtain
\begin{multline*}
C_{\overline{M,N}}  = \frac{1}{\phi^{*}(q)}  \sum_{d | q} \phi(d) \mu(\frac{q}{d}) \sum_{l} \frac{1}{l^{1 - \gamma + \delta}} \sumstar_{h \shortmod{l}}
\sum_{f \equiv 0 \shortmod{d}} 
\sum_{m}  \frac{\sigma_{\alpha - \beta}(m) e\left(\frac{mh}{l} \right)e\left(\frac{-hf}{l} \right)}{m^{\half + \alpha}}
\\
\left(\frac{1}{2 \pi i}\right)^4 \int_{(c_w)} \int_{(c_{u_1})} \int_{(c_{u_2})} \int_{(c_s)} 
\frac{q^{2s} M^{u_1} N^{u_2} H_1(s, u_1, u_2, w) }{l^w m^{s +u_1 } (f -m)^{\half + \gamma +s + u_2 -\frac{w}{2}} }  ds du_1 du_2 dw,
\end{multline*}
where $H_1$ is given by \eqref{eq:Hthree}.

To separate the variables $m$ and $f$ in $(f-m)^{-1/2 - \gamma - s -u_2 + \frac{w}{2}}$ we use the following formula:
\begin{align*}
\frac{1}{2 \pi i} \int_{(c_v)} \frac{\Gamma(v) \Gamma(1-b)}{ \Gamma(1-b+v)} x^{-v} dv 
=\begin{cases}
(1 - x)^{-b}, \quad &0  < x < 1 \\
0, \quad &1 < x.
\end{cases}
,
\end{align*}
valid for $0 < c_v$, $\text{Re}(b) <1$, from 17.43.22 of \cite{GR}.  The integration converges absolutely provided $\text{Re}(b) < 0$ using Stirling's approximation, namely
\begin{equation*}
\frac{\Gamma(v) \Gamma(1-b)}{\Gamma(1-b + v)} \ll (1+ |\text{Im}(v)|)^{\text{Re}(b) - 1},
\end{equation*}
where the implied constant depends polynomially on $b$ (the dependence on $b$ could be explicitly stated but it turns out not to be relevant in this work).
Thus we have
\begin{multline*}
C_{\overline{M,N}} = \frac{1}{\phi^{*}(q)}  \sum_{d | q} \phi(d) \mu(\frac{q}{d}) \sum_{l} \frac{1}{l^{1 - \gamma + \delta}}
\sumstar_{h \shortmod{l}}    
\sum_{\substack{f \equiv 0 \shortmod{d}}}{} \frac{e\left(\frac{-hf}{l} \right) }{f^{\half + \gamma}}
\sum_{m}  \frac{\sigma_{\alpha - \beta}(m) e\left(\frac{mh}{l}\right)}{m^{\half + \alpha}} 
\\
\left(\frac{1}{2 \pi i} \right)^5 \int_{(c_s)} \int_{(c_w)} \int_{(c_{u_2})} \int_{(c_{u_1})} \int_{(c_v)}   \frac{q^{2s} M^{u_1} N^{u_2} \overline{H}(s, u_1, u_2, v, w)}{f^{ s+ u_2 - v - \frac{w}{2}} l^w m^{s + u_1 + v}  }dv du_1 du_2 dw ds,
\end{multline*}
where
\begin{equation*}
\overline{H}(s,u_1,u_2,v,w) = \frac{\Gamma(v) \Gamma(\half - \gamma - s  - u_2 +\frac{w}{2})}{\Gamma(\half - \gamma - s  - u_2 +v +\frac{w}{2})} H_1(s,u_1,u_2,w).
\end{equation*}
Note that this expression is identical to \eqref{eq:blackvitamin} except that the ratio of gamma factors arising from the separation of variables is different.  We essentially follow the same outline used to compute $C_{M,N}$ but some arguments must be altered due to changes in location of the poles of these gamma functions.  We initially take $c_s = \half$, $c_{u_1} = c_{u_2} = 0$, $c_v = \varepsilon$, and $c_w = 2 + \varepsilon$ so that all the integrals and the sum over $m$ converge absolutely.  Then we write the sum over $m$ in terms of the Estermann $D$-function and move $c_s$ to $\varepsilon$, thereby crossing the two poles of $D$.  Let $E_{\overline{M,N}}$ be the integral along the new line of integration, and $P_{\overline{M,N}}(\alpha,\beta,\gamma,\delta)$ be the contribution of the pole at $\half + s + u_1 + v = 1$.  We compute $E_{\overline{M,N}}$ in Section \ref{section:errorterms} and proceed with the computation of the main terms.

\subsection{The main terms}
The same argument used in the proof of Lemma \ref{lemma:POmainterm} shows how we may execute the summation over $M$ and $N$.  Furthermore, with the same method of proof as Lemma \ref{lemma:claimedupperboundonP} we have
\begin{equation}
P_{\overline{M,N}} \ll M^{\half} N^{-\half} q^{\varepsilon}.
\end{equation}
Let $\overline{P}(\alpha,\beta,\gamma,\delta) := \sum_{M,N} P_{\overline{M,N}}(\alpha,\beta,\gamma,\delta)$.
The analog of Lemma \ref{lemma:POmainterm} is
\begin{mylemma}
We have
\begin{multline}
\label{eq:mt3}
\overline{P}(\alpha,\beta,\gamma,\delta) = \zeta(1 - \alpha + \beta)     
\Big(\frac{1}{2 \pi i} \Big)^2 
\int_{(\varepsilon)} \int_{(\varepsilon)} 
 q^{-\alpha - \gamma + \frac{w}{2}}
\frac{G(v)G(w) g(v)}{ vw}
\zeta(1 - \gamma + \delta + w)
\\
\frac{\Gamma(\half -\alpha  -v)\Gamma(\half - \gamma -v + \frac{w}{2})}{\Gamma(1 - \alpha - \gamma -2 v+ \frac{w}{2})}
\frac{\zeta( \alpha + \gamma +2v - \frac{w}{2}) \zeta(1 + \beta + \delta +2v + \frac{w}{2})}{\zeta(2-\alpha + \beta -\gamma + \delta + w)}
dv dw
\\
+ O(q^{-1 +\varepsilon}).
\end{multline}
\end{mylemma}
Note that the pole of $\zeta(\alpha + \gamma + 2v - \frac{w}{2})$ is cancelled by one of the gamma factors.
\begin{proof}
To begin, $\overline{P}$ is given by the residue of the following integral at $s + v = \half - \alpha$:
\begin{multline*}
\frac{1}{\phi^{*}(q)}  \sum_{d | q} \phi(d) \mu(\frac{q}{d}) \sum_{l} \frac{1}{l^{1 - \gamma + \delta}}
\sumstar_{h \shortmod{l}}    
\sum_{\substack{f \equiv 0 \shortmod{d}}}{} \frac{e\left(\frac{-hf}{l} \right) }{f^{\half + \gamma}}
\\
\left(\frac{1}{2 \pi i} \right)^3 \int_{(c_s)} \int_{(c_w)} \int_{(c_v)}   \frac{q^{2s} }{f^{s - v - \frac{w}{2} } l^w  }
D(\tfrac12 + \alpha + s  +v, \alpha - \beta, \frac{h}{l})
\overline{H}_3(s, v, w)
dv  dw ds,
\end{multline*}
where
\begin{equation*}
\overline{H}_3(s,v,w) = \frac{\Gamma(v) \Gamma(\half - \gamma - s + \frac{w}{2})}{\Gamma(\half - \gamma - s +v + \frac{w}{2})} 
\frac{G(s) G(w)}{s w} g(s) \zeta(1 - \gamma + \delta + w).
\end{equation*}

Now we compute $\overline{P}(\alpha,\beta,\gamma,\delta)$.
Computing the residue and changing variables $v \rightarrow \half - \alpha - v$ gives the following analog of \eqref{eq:spoon}
\begin{multline*}
\overline{P} = \frac{\zeta(1 - \alpha + \beta)}{\phi^{*}(q)} \sum_{d | q} \phi(d) \mu(q/d) \sum_{l} \frac{1}{l^{2 - \alpha + \beta - \gamma + \delta}}
\sumstar_{h \shortmod{l}}   
\sum_{\substack{f \equiv 0 \shortmod{d}}}{} e\left(\frac{-hf}{l} \right)
\\
\left(\frac{1}{2 \pi i} \right)^2  \int_{(c_w)} \int_{(c_v)}   \frac{q^{2 v} }{f^{ \alpha + \gamma  + 2v - \frac{w}{2} } l^w  }
\overline{H}_3(v, \tfrac12 - \alpha - v, w)
dv dw,
\end{multline*}
where note the ratio of gamma factors in $\overline{H}_3$ is
\begin{equation*}
\frac{\Gamma(\half - \alpha - v) \Gamma( \half - \gamma -v + \frac{w}{2})}{\Gamma(1 - \alpha - \gamma -2v + \frac{w}{2})}.
\end{equation*}
Now we take $c_w = 3 \varepsilon$ and move $c_v$ to $\half + \varepsilon$.  The pole of $\Gamma(\half - \alpha -v)$ is cancelled by a zero of $G$, and the sum over $f$ converges absolutely, so we may borrow the same computations used in the proof of Lemma \ref{lemma:POmainterm} to get \eqref{eq:mt3}.
\end{proof}

\section{Assembling the main terms}
\label{section:assemble}
We now begin to assemble the various main terms to form a nicer expression.  We briefly recall how the main terms decompose.  Recall \eqref{eq:MAFE}, \eqref{eq:switch}, and \eqref{eq:A1q}.  In general, \eqref{eq:switch} shows how to derive results for the ``second part'' of the approximate functional equation.  Recall the combinatorial dissection \eqref{eq:decomposition}, as well as the computation of $A_D$ from Lemma \ref{lem:YandJ}.  Then we wrote $B_O = \sum_{M,N} B_{M,N}$ and recall \eqref{eq:BC}, 
\eqref{eq:CPE}, and Lemma \ref{lemma:POmainterm}.  This gives $B_O$ as the sum of four main terms, and we similarly express $B_O'(\alpha, \beta, \gamma, \delta) = B_O(\gamma, \delta, \alpha, \beta)$.  The computations are similar for the dual terms as in Section \ref{section:dual}.

\subsection{Combining $A_O$ and $A_{\overline{O}}$}
To begin, let $Q=Q(\alpha,\beta,\gamma,\delta)$ be the sum of \eqref{eq:POmainterm}, \eqref{eq:mt2}, and \eqref{eq:mt3}.  
From the above discussion, the total contribution of main terms from $B_O$, $B_O'$, and $A_{\overline{O}}$ is then 
\begin{equation}
\label{eq:vacuumcleaner}
Q(\alpha, \beta, \gamma, \delta) + Q(\beta, \alpha, \gamma, \delta) + Q(\alpha, \beta, \delta, \gamma) + Q(\beta, \alpha, \delta, \gamma).
\end{equation}
This grouping of the main terms is suggested by unpublished work of Hughes on the Riemann zeta function \cite{H}.
\begin{mylemma}
\label{lemma:half}
We have
\begin{multline}
\label{eq:half}
Q(\alpha,\beta, \gamma, \delta) =  \left[ \frac{\zeta(1 - \alpha + \beta)\zeta(1 - \gamma + \delta)}{\zeta(2 -\alpha + \beta - \gamma + \delta)}
 \left(\frac{q}{\pi}\right)^{-\alpha - \gamma}
 \frac{1}{2 \pi i} \int_{(\frac14)} 
\frac{G(s)}{s} \pi^{2s} g_{\alpha, \beta, \gamma, \delta}(s)
\right.
\\
\left. \zeta(1 - \alpha - \gamma - 2s) \zeta(1 + \beta + \delta +2s)
\frac{\Gamma\left(\frac{\half - \alpha - s}{2}\right)}{\Gamma\left(\frac{\half + \alpha + s}{2}\right)}
\frac{\Gamma\left(\frac{\half - \gamma - s}{2}\right)}{\Gamma\left(\frac{\half + \gamma + s}{2}\right)} 
ds \right] + O(q^{-1/3+\varepsilon}).
\end{multline}
\end{mylemma}
\begin{proof}
By its definition,
\begin{multline}
\label{eq:Q}
Q =   
\left(\frac{1}{2 \pi i} \right)^2 
\int_{(\frac14)} \int_{(\varepsilon)} 
 q^{-\alpha - \gamma + \frac{w}{2}}  
\frac{G(v)G(w) g(v)}{vw} R(v,w)
\\ 
 \frac{\zeta( \alpha + \gamma +2v - \frac{w}{2}) \zeta(1 + \beta + \delta +2v + \frac{w}{2})}{\zeta(2-\alpha + \beta -\gamma + \delta + w)}
dw dv,
\end{multline}
where
\begin{multline*}
R(v,w) = \zeta(1-\alpha + \beta) \zeta(1 - \gamma + \delta + w) \frac{\Gamma(\half -\alpha -v) \Gamma(\alpha + \gamma + 2v - \frac{w}{2})}{\Gamma(\half + \gamma + v - \frac{w}{2})}
\\
+ \zeta(1-\alpha + \beta + w) \zeta(1 - \gamma + \delta) \frac{\Gamma(\half -\gamma -v) \Gamma(\alpha + \gamma + 2v - \frac{w}{2})}{\Gamma(\half + \alpha + v - \frac{w}{2})}
\\
+\zeta(1-\alpha + \beta) \zeta(1 - \gamma + \delta + w) \frac{\Gamma(\half -\alpha -v) \Gamma(\half - \gamma -v + \frac{w}{2})}{\Gamma(1 -\alpha - \gamma - 2v + \frac{w}{2})}.
\end{multline*}
We moved the line of integration $c_v$ to $\frac{1}{4}$ without encountering any poles.
Now move the line of integration over $w$ to $-1 +  \varepsilon$, crossing poles at 
$w=0$, and $\frac{w}{2} = -\half + \gamma + v$ only (recall Definition \ref{def:G}).
The contribution from the pole at $w = -1 + 2\gamma + 2v$ is
\begin{multline}
\label{eq:weird}
\frac{\phi(q)}{\phi^{*}(q)}    
\frac{1}{2 \pi i} 
\int_{(\frac14)} 
 q^{-\half - \alpha +v}  
\frac{G(v)G(-1 + 2\gamma + 2v) g(v)}{v(-1 + 2\gamma + 2v)} 
\zeta(1-\alpha + \beta) \zeta( \gamma + \delta + 2v) 
\\
\prod_{p | q} \left(1 - p^{-\half + \alpha + v }\right) \frac{\zeta(\half + \alpha + v) \zeta(\half + \beta + \gamma+ \delta +3v)}{\zeta(1-\alpha + \beta +\gamma + \delta +2v)}
dw dv.
\end{multline}
It is clear that this integral is bounded by $O(q^{-1/3 + \varepsilon})$, as can be seen by moving $v$ to $\varepsilon$ and capturing the pole at $3v \approx \half$.

The pole at $w=0$ gives
\begin{equation*} 
\frac{1}{2 \pi i} 
\int_{(\frac14)}
 q^{-\alpha - \gamma}
\frac{G(v) g(v)}{v} R(v, 0)
\frac{\zeta( \alpha + \gamma +2v) \zeta(1 + \beta + \delta +2v)}{\zeta(2-\alpha + \beta -\gamma + \delta)}
 dv.
\end{equation*}
Now apply the functional equation to $\zeta(\alpha + \gamma + 2v)$ to express this as
\begin{equation*}
\frac{1}{2 \pi i} 
\int_{(\frac14)}
\left(\frac{q}{\pi}\right)^{-\alpha - \gamma}  
\frac{G(v) g(v) R(v,0)}{\pi^{\half-2v} v} 
  \frac{\Gamma(\frac{1 - \alpha - \gamma - 2v}{2})}{\Gamma(\frac{ \alpha + \gamma + 2v}{2})}
\frac{\zeta(1- \alpha - \gamma -2v ) \zeta(1 + \beta + \delta +2v)}{\zeta(2-\alpha + \beta -\gamma + \delta)}
 dv.
\end{equation*}
Note
\begin{equation*}
R(v, 0) = \zeta(1-\alpha + \beta) \zeta(1 - \gamma + \delta) \Gamma_{\alpha, \gamma}(v), 
\end{equation*}
where
\begin{multline}
\label{eq:butter}
\Gamma_{\alpha, \gamma}(v) = 
\frac{\Gamma(\half -\alpha -v) \Gamma(\alpha + \gamma + 2v)}{\Gamma(\half + \gamma + v)}
+ \frac{\Gamma(\half -\gamma -v) \Gamma(\alpha + \gamma + 2v )}{\Gamma(\half + \alpha + v)}
\\
+ \frac{\Gamma(\half -\alpha -v) \Gamma(\half - \gamma -v )}{\Gamma(1 -\alpha - \gamma - 2v)}.
\end{multline}
To finish the proof of Lemma \ref{lemma:half} we use the crucial identity
\begin{equation}
\Gamma_{\alpha, \gamma}(v) = 
\pi^{\half} 
\frac{\Gamma\left(\frac{ \alpha + \gamma +2v}{2}\right) }{\Gamma\left(\frac{1 - \alpha - \gamma -2v}{2}\right)}
\frac{\Gamma\left(\frac{\half - \alpha - v}{2}\right)}{\Gamma\left(\frac{\half + \alpha + v}{2}\right)}
\frac{\Gamma\left(\frac{\half - \gamma - v}{2}\right)}{\Gamma\left(\frac{\half + \gamma + v}{2}\right)},
\end{equation}
which is deduced from the following Lemma with $a= \half - \alpha - v$ and $b=\half - \gamma -v$.
\end{proof}

\begin{mylemma}
\label{lemma:gamma}
For any $a, b \in \mc$, $a,b, a+b \not \in \mz$,
\begin{equation}
\label{eq:gamma}
\frac{\Gamma(a ) \Gamma(1-a-b)}{ \Gamma(1-b)} 
+ \frac{\Gamma(b) \Gamma(1-a-b)}{ \Gamma(1-a)} 
+ \frac{\Gamma(a) \Gamma(b)}{ \Gamma(a+b)}
=
\pi^{\half} 
\frac{\Gamma\left(\frac{ 1-a-b}{2}\right) }{\Gamma\left(\frac{a+b}{2}\right)}
\frac{\Gamma\left(\frac{a}{2}\right)}{\Gamma\left(\frac{1-a}{2}\right)}
\frac{\Gamma\left(\frac{b}{2}\right)}{\Gamma\left(\frac{1-b}{2}\right)}.
\end{equation}
\end{mylemma}
\begin{proof}[Proof of Lemma \ref{lemma:gamma}]
Using
\begin{equation*}
\frac{\Gamma(\frac{s}{2})}{\Gamma(\frac{1-s}{2})} = \pi^{-\half} 2^{1-s} \cos({\textstyle \frac{\pi s}{2}}) \Gamma(s),
\end{equation*}
the right hand side of \eqref{eq:gamma} equals
\begin{equation}
\label{eq:LHStrig}
2 \frac{\cos(\frac{\pi a}{2}) \cos(\frac{\pi b}{2}) \Gamma(a) \Gamma(b)}{\cos(\frac{\pi (a + b)}{2}) \Gamma(a+b)}.
\end{equation}
Using a series of standard gamma function and trigonometric identities, the left hand side is
\begin{gather*}
\frac{\Gamma(a)\Gamma(b)}{\Gamma(a + b)} 
\left(
	\frac{\Gamma(a+b) \Gamma(1-a-b)}{\Gamma(b) \Gamma(1-b)} 
	+ \frac{\Gamma(a+b) \Gamma(1-a-b)}{\Gamma(a) \Gamma(1-a)}
	+ 1
\right)
\\
= \frac{\Gamma(a)\Gamma(b)}{\Gamma(a + b)} 
\left(\frac{\sin(\pi a) + \sin( \pi b) + \sin(\pi(a + b))}{\sin(\pi(a+b))} \right)
\\
= 2 \frac{\Gamma(a)\Gamma(b)}{\Gamma(a + b)} \frac{\sin(\frac{\pi}{2}(a+b)) \cos(\frac{\pi}{2}(a-b)) + \sin(\frac{\pi}{2}(a+b)) \cos(\frac{\pi}{2}(a+b)) }{\sin(\pi(a+b))}
\\
= 4 \frac{\Gamma(a)\Gamma(b)}{\Gamma(a + b)} \frac{\sin(\frac{\pi}{2}(a+b)) \cos(\frac{\pi a}{2}) \cos(\frac{\pi b}{2})  }{\sin(\pi(a+b))},
\end{gather*}
which simplifies to give \eqref{eq:LHStrig} from the sine double angle formula.
\end{proof}

\subsection{Combining $Q$ and $Q_{-}$}
In this section we calculate the terms corresponding to $Q$ that create the main terms of $M_{-1}$.  Recall that these are obtained from $Q$ by switching the signs of $\alpha, \beta, \gamma, \delta$ and multiplying by $X_{\alpha,\beta,\gamma,\delta}$.
We combine $Q(\alpha, \beta, \gamma, \delta)$ and $Q_{-}(\beta, \alpha, \delta, \gamma)$.
\begin{mylemma}
\label{lemma:QandQminus}
We have
\begin{equation}
\label{eq:QandQminus}
Q(\alpha, \beta, \gamma, \delta) + Q_{-}(\beta, \alpha, \delta, \gamma)
=
U(\alpha,\beta,\gamma,\delta)
+ O(q^{-1/3 + \varepsilon}),
\end{equation}
where
\begin{equation}
U(\alpha,\beta,\gamma,\delta) = 
X_{\alpha,\gamma} \frac{\zeta(1 - \alpha + \beta) \zeta(1 - \alpha - \gamma ) \zeta(1 + \beta + \delta)\zeta(1 - \gamma + \delta)}{\zeta(2 -\alpha + \beta - \gamma + \delta)}.
 \end{equation}
\end{mylemma}
Note that $U(\alpha,\beta,\gamma,\delta)$ is one of the middle four terms in Conjecture \ref{conjecture:five}. Hence, adding \eqref{eq:QandQminus} according to \eqref{eq:vacuumcleaner} to the contribution of the diagonal terms given by Lemma \ref{lem:YandJ} will form the quantity on the right hand side of Conjecture \ref{conjecture:five}.

\begin{proof}
It is a matter of bookkeeping to modify \eqref{eq:half} to see  
\begin{multline}
\label{eq:Q-}
Q_{-}(\beta, \alpha, \delta, \gamma) =  X_{\alpha,\beta,\gamma,\delta}  
\frac{\zeta(1 - \alpha + \beta)\zeta(1 - \gamma + \delta)}{\zeta(2 -\alpha + \beta - \gamma + \delta)}
 \left(\frac{q}{\pi}\right)^{\beta + \delta} \frac{1}{2 \pi i} \int_{(\frac14)} \frac{G(s)}{s} \pi^{2s}
 \\
  \zeta(1 - \alpha - \gamma +2s) \zeta(1 + \beta + \delta - 2s) 
\frac{\Gamma\left(\frac{\half + \beta - s}{2}\right)}{\Gamma\left(\frac{\half - \beta + s}{2}\right)}
\frac{\Gamma\left(\frac{\half - \delta + s}{2}\right)}{\Gamma\left(\frac{\half - \delta + s}{2}\right)}
  g_{-\alpha, -\beta, -\gamma, -\delta}(s) 
ds
+ O(q^{-1/3+\varepsilon}).
\end{multline}
Let $I$ be the main term in \eqref{eq:half}, and $I_{-}$ be the main term in \eqref{eq:Q-}.  We work with $I$ by moving the line of integration to $-\frac14$, passing a pole at $s=0$ giving a residue which we easily compute to be $U(\alpha,\beta,\gamma,\delta)$.
Write $I = U(\alpha, \beta, \gamma, \delta) + I'$, where $I'$ is the new integral.  We now show $I' = -I_{-}$, which will complete the proof.

Apply the change of variables $s \rightarrow -s$ to give
\begin{multline*}
I' = -
 \frac{\zeta(1 - \alpha + \beta)\zeta(1 - \gamma + \delta)}{\zeta(2 -\alpha + \beta - \gamma + \delta)}
 \left(\frac{q}{\pi}\right)^{-\alpha - \gamma}
 \\
 \frac{1}{2 \pi i} \int_{(\frac14)} \zeta(1 - \alpha - \gamma + 2s) \zeta(1 + \beta + \delta -2s)
\frac{\Gamma\left(\frac{\half - \alpha + s}{2}\right)}{\Gamma\left(\frac{\half + \alpha - s}{2}\right)}
\frac{\Gamma\left(\frac{\half - \gamma + s}{2}\right)}{\Gamma\left(\frac{\half + \gamma - s}{2}\right)}
 \frac{G(s)}{s} \pi^{-2s} g_{\alpha, \beta, \gamma, \delta}(-s) 
ds.
\end{multline*} 
Now we claim that
\begin{multline*}
\left(\frac{q}{\pi}\right)^{-\alpha - \gamma} \frac{\Gamma\left(\frac{\half - \alpha + s}{2}\right)}{\Gamma\left(\frac{\half + \alpha - s}{2}\right)}
\frac{\Gamma\left(\frac{\half - \gamma + s}{2}\right)}{\Gamma\left(\frac{\half + \gamma - s}{2}\right)}  \pi^{-2s} g_{\alpha, \beta, \gamma, \delta}(-s)
\\
=
X_{\alpha,\beta,\gamma,\delta} \left(\frac{q}{\pi}\right)^{\beta + \delta}
\frac{\Gamma\left(\frac{\half + \beta - s}{2}\right)}{\Gamma\left(\frac{\half - \beta + s}{2}\right)}
\frac{\Gamma\left(\frac{\half + \delta + s}{2}\right)}{\Gamma\left(\frac{\half - \delta + s}{2}\right)}
\pi^{2s} g_{-\alpha, -\beta, -\gamma, -\delta}(s),
\end{multline*}
which implies $I' = -I_{-}$.
Each of these terms is a product of terms each depending on exactly one of $\alpha$, $\beta$, $\gamma$, or $\delta$, as well as a certain power of $\pi$ not depending on the shifts, so it suffices to check this identity at each such factor.  The cases $\alpha$ and $\gamma$ are the same and so are $\beta$ and $\delta$.  The case of $\alpha$ follows from
\begin{equation*}
\left(\frac{q}{\pi}\right)^{-\alpha} \frac{\Gamma\left(\frac{\half - \alpha + s}{2}\right)}{\Gamma\left(\frac{\half + \alpha - s}{2}\right)} \frac{\Gamma\left(\frac{\half + \alpha - s}{2}\right)}{\Gamma\left(\frac{\half + \alpha}{2}\right)}
=
X(\tfrac12 + \alpha)  \frac{\Gamma\left(\frac{\half - \alpha + s}{2}\right)}{\Gamma\left(\frac{\half - \alpha}{2}\right)},
\end{equation*}
and the analogous formula for $\beta$ is
\begin{equation*}
\frac{\Gamma\left(\frac{\half + \beta - s}{2}\right)}{\Gamma\left(\frac{\half + \beta}{2}\right)}
=
X(\tfrac12 + \beta) \left(\frac{q}{\pi}\right)^{\beta}
\frac{\Gamma\left(\frac{\half + \beta - s}{2}\right)}{\Gamma\left(\frac{\half - \beta + s}{2}\right)}
\frac{\Gamma\left(\frac{\half - \beta + s}{2}\right)}{\Gamma\left(\frac{\half - \beta}{2}\right)}.
\end{equation*}
We also check that the power of $\pi$ is the same on both sides.
\end{proof}

\subsection{A note on odd characters}
\label{section:odd}
In this work we concentrate almost exclusively on the even characters in the proof of Theorem \ref{thm:mainresult}.  The contribution of the odd characters carries through in the same way with slight changes.  The only differences between the odd and even characters is that the $X$ factors are different, and for the odd characters, the `dual' terms are subtracted rather than added.  The estimations of the error terms carry through as before; the only differences arise in the calculation of the main terms.  Some thought shows that the evaluation of the main terms for the odd characters diverges from that of the even characters starting with the analog of \eqref{eq:butter}; for the odd character case the third term is subtracted rather than added.  The only essential difference is the use of the following gamma identity instead of Lemma \ref{lemma:gamma}.
\begin{mylemma}
\label{lemma:gamma2}
For any $a, b \in \mc$, $a,b, a+b \not \in \mz$,
\begin{multline}
\frac{\Gamma(a ) \Gamma(1-a-b)}{ \Gamma(1-b)} 
+ \frac{\Gamma(b) \Gamma(1-a-b)}{ \Gamma(1-a)} 
- \frac{\Gamma(a) \Gamma(b)}{ \Gamma(a+b)}
=
\pi^{\half} 
\frac{\Gamma\left(\frac{ 1-a-b}{2}\right) }{\Gamma\left(\frac{a+b}{2}\right)}
\frac{\Gamma\left(\frac{1+a}{2}\right)}{\Gamma\left(\frac{2-a}{2}\right)}
\frac{\Gamma\left(\frac{1+b}{2}\right)}{\Gamma\left(\frac{2-b}{2}\right)}.
\end{multline}
\end{mylemma}
The proof is similar to that of Lemma \ref{lemma:gamma} so we omit the details.  

\section{Treating the error terms}
\label{section:errorterms}
In this section we prove Theorem \ref{thm:Ebound}.  Since the forthcoming estimations become techncial, it may be helpful to know that the quality of the error term in Theorem \ref{thm:Ebound} can be predicted by careful scrutiny of the formulas in Theorems 3 and 4 of \cite{Motohashi1}.  

In order to clean the upcoming formulas, we set all the shift parameters $\alpha, \beta$, etc., equal to $0$.  All the arguments can easily be generalized to handle sufficiently small non-zero values without degrading the results.  
Recall that $E_{M,N}$ is given by the right hand side of \eqref{eq:POintegral}
with contours of integration defined by
\begin{equation}
\label{eq:lines}
c_s = c_w = c_v = \varepsilon, \quad c_{u_1} = c_{u_2} = 0.
\end{equation}
\subsection{Reduction to Kloosterman sums}
\label{section:reduction}
Apply the functional equation \eqref{eq:FE} to $D$ to obtain $E_{M,N} = E_{+} + E_{-}$, where
\begin{multline}
\label{eq:Epm}
E_{\pm} = \frac{1}{\phi^{*}(q)}  \sum_{d | q} \phi(d) \mu(\frac{q}{d}) \sum_{l} \frac{1}{l}
\sumstar_{h \shortmod{l}}    
\sum_{\substack{f \equiv 0 \shortmod{d}}}{} \frac{e\left(\frac{hf}{l} \right) }{f^{\half}}
\\
\left(\frac{1}{2 \pi i} \right)^5 \int_{(c_s)} \int_{(c_w)} \int_{(c_{u_2})} \int_{(c_{u_1})} \int_{(c_v)}   \frac{q^{2s} M^{u_1} N^{u_2} }{f^{s+ u_2 - v - \frac{w}{2} } l^{2s + 2u_1 + 2v +w}  }
\\
D(\tfrac12 - s - u_1 -v, 0, \frac{\pm \overline{h}}{l})
H_{\pm}(s, u_1, u_2, v, w)
dv du_1 du_2 dw ds,
\end{multline}
and where
\begin{equation}
\label{eq:Hpmdefinition}
H_{\pm}(s,u_1, u_2 ,v,w) = 2 (2\pi)^{-1 + 2s + 2u_1 + 2v} \Gamma(\tfrac12 -s -u_1 - v)^2 H(s, u_1, u_2, v, w)
 S_{\pm},
\end{equation}
with
\begin{equation}
\label{eq:Splusminus}
S_{+} = 1, \quad S_{-} = \sin(\pi(s + u_1 +v)).
\end{equation}
Here $H_{\pm}$ has rapid decay in all variables since $H$ does and since the exponential decay of the gamma factors cancels the exponential growth of $S_{-}$.

Now move $c_{u_1}$ to $-1$ and expand $D$ into absolutely convergent Dirichlet series and execute the sum over $h$ in terms of Kloosterman sums to get
\begin{multline}
\label{eq:theonion}
E_{\pm} = \frac{1}{\phi^{*}(q)}  \sum_{d | q} \phi(d) \mu(\frac{q}{d}) \sum_{l} \frac{1}{l}  
\sum_{\substack{f \equiv 0 \shortmod{d}}}{} \sum_m \frac{S(m,\pm f;l) }{f^{\half} m^{\half}}
\left(\frac{1}{2 \pi i} \right)^5 \int_{(c_s)} \int_{(c_w)} 
\\
\int_{(c_{u_2})} \int_{(c_{u_1})} \int_{(c_v)}   \frac{q^{2s} M^{u_1} N^{u_2} m^{s + u_1 + v}}{f^{s+ u_2 - v - \frac{w}{2} } l^{2s + 2u_1 + 2v +w}  }
H_{\pm}(s, u_1, u_2, v, w)
dv du_1 du_2 dw ds.
\end{multline}
As an aside, we mention that using only the Weil bound for Kloosterman sums, we obtain the bound
\begin{equation*}
E_{M,N} \ll q^{-1 + \varepsilon} M^{-\half} N^{\frac54},
\end{equation*}
which can be seen easily by taking 
\begin{equation*}
c_s = \varepsilon, \quad c_v = \varepsilon, \quad c_{u_1} = -\tfrac12 - 3\varepsilon, \quad c_{u_2} = \tfrac{5}{4} + 3\varepsilon, \quad c_w = \tfrac32 + 3\varepsilon.
\end{equation*}
Using $MN \ll q^{2 + \varepsilon}$ shows
\begin{equation*}
E_{M,N} \ll M^{-1} N^{\frac34} q^{\varepsilon},
\end{equation*}
which is nontrivial only for $N \ll q^{\frac87 - \varepsilon}$, $M \gg q^{\frac67 + \varepsilon}$.  This is insufficient to combine with Theorem \ref{prop:faraway} to cover all ranges; for example, $M=q^{2/3}$, $N= q^{4/3}$ would not be covered.

To do better we shall obtain additional savings coming from cancellation in the sum of Kloosterman sums by the use of the Kuznetsov formula.

\subsection{Preparation for application of Kuznetsov}
Let
\begin{multline*}
r_{\pm}(x) =
\left(\frac{1}{2 \pi i} \right)^5 \int_{(c_s)} \int_{(c_w)} \int_{(c_{u_2})} \int_{(c_{u_1})} \int_{(c_v)}   \frac{q^{2s} M^{u_1} N^{u_2} }{f^{2s + u_1 + u_2 } m^{\frac{w}{2}} }
\left(\frac{x}{4\pi}\right)^{2s + 2u_1 +2v + w}
\\
H_{\pm}(s, u_1, u_2, v, w)
dv du_1 du_2 dw ds,
\end{multline*}
with contours of integration given by \eqref{eq:lines}.
Then
\begin{equation*}
E_{\pm} = \frac{1}{\phi^{*}(q)}  \sum_{d | q} \phi(d) \mu(\frac{q}{d})    
\sum_{\substack{f \equiv 0 \shortmod{d}}}{} \sum_m \frac{d(m)}{ f^{\half} m^{\half}} 
\sum_{l} \frac{S(m, \pm f;l)}{l} r_{\pm} \left(\frac{ 4 \pi \sqrt{mf}}{l} \right).
\end{equation*}
By taking $c_s = \half - 2\varepsilon$, we get $r_{\pm}(x) \ll x^{1 - \varepsilon}$.  By taking $c_{u_1} = -A$  with $A$ large, we see that $r^{(j)}(x) \ll_{j,B} (1 + x)^{-B}$ for $j= 0,1,2,\ldots$ and any $B > 0$.  These conditions are sufficient for the application of the Kuznetsov formula, Theorem \ref{thm:kuznetsov}.

The next step is to apply the Kuznetsov formula to the sum over $l$.  We write $E_{\pm} = E_{m \pm } + E_{c \pm} + E_{h \pm}$ to correspond to the Maass forms, the Eisenstein series, and the holomorphic forms (of course $E_{h-} = 0$).  We show how to estimate $E_{m\pm}$ and $E_{c\pm}$, since $E_{h+}$ is smaller and easier to handle (for example, see Section 5 of \cite{Motohashi1}).

\subsection{Integral transforms}
At this point we manipulate the various integral transforms of $r_{\pm}$ that we require for the application of the Kuznetsov formula.  We require $M_{r_{+}}(t)$, $K_{r-}(t)$, and $N_{r_{+}}(k)$, in the notation of \cite{IK}, Theorems 16.5 and 16.6.
We have
\begin{multline*}
M_{r+}(t) =  \left(\frac{1}{2 \pi i} \right)^5 \int_{(c_s)} \int_{(c_w)} \int_{(c_{u_2})} \int_{(c_{u_1})} \int_{(c_v)}   \frac{q^{2s}  M^{u_1} N^{u_2} H_{+}(s, u_1, u_2, v, w)}{(4\pi)^{2s +2u_1 + 2v + w} f^{2s +u_1 + u_2} m^{\frac{w}{2}}}  
\\
\frac{\pi i}{\sinh{2 \pi t}} \int_0^{\infty} \left(J_{2it}(x) - J_{-2it}(x) \right)
x^{2s +2u_1 +  2v +w} \frac{dx}{x}
dv du_1 du_2 dw ds.
\end{multline*}
Now use the formula
\begin{equation*}
\int_0^{\infty} J_{\nu}(x) x^s \frac{dx}{x} = 2^{s-1} \frac{\Gamma(\frac{s + \nu}{2})}{\Gamma(\frac{\nu-s +2}{2})},
\end{equation*}
valid for $\frac{3}{2} > \text{Re}(s) > -\text{Re}(\nu)$
(see 6.561.14) of \cite{GR}) to see that
\begin{equation*}
\int_0^{\infty} (J_{2it}(x) - J_{-2it}(x)) x^{\lambda} \frac{dx}{x} = 2^{\lambda -1} \left(\frac{\Gamma(\frac{\lambda}{2}  + it)}{\Gamma(1 -\frac{\lambda}{2}  + it))} - \frac{\Gamma(\frac{\lambda}{2}  - it))}{\Gamma(1 -\frac{\lambda}{2}  - it)} \right).
\end{equation*}
We simplify this expression using
\begin{mylemma}
\label{lemma:besselgamma}
\begin{equation}
\frac{\Gamma(a + ir)}{\Gamma(1-a + ir)} - \frac{\Gamma(a - ir)}{\Gamma(1-a - ir)} = -\frac{2i}{\pi} \sinh(\pi r) \cos(\pi a) \Gamma(a + ir) \Gamma(a - ir).
\end{equation}
\end{mylemma}
We omit the proof since it is easy and standard.
Thus
\begin{equation*}
\frac{\pi i 2^{\lambda -1}}{\sinh{2 \pi t}}   \left(\frac{\Gamma(\frac{\lambda}{2}  + it)}{\Gamma(1 -\frac{\lambda}{2}  + it))} - \frac{\Gamma(\frac{\lambda}{2}  - it))}{\Gamma(1 -\frac{\lambda}{2}  - it)} \right) 
=  \frac{2^{\lambda-1}}{\cosh{\pi t}}\cos(\frac{ \pi \lambda}{2}) \Gamma(\frac{\lambda}{2} + it) \Gamma(\frac{\lambda}{2}  - it).
\end{equation*}
Hence
\begin{equation*}
M_{r+}(t) =\left(\frac{1}{2 \pi i} \right)^5 
\int \dots \int  \frac{q^{2s} M^{u_1} N^{u_2}}{f^{2s +u_1 + u_2} m^{\frac{w}{2}}}   \frac{\widehat{H}_{+}(s, u_1, u_2, v, w;t)}{\cosh(\pi t)}
dv du_1 du_2 dw ds,
\end{equation*}
where
\begin{multline*}
\widehat{H}_{+}(s,u_1, u_2,v,w; t) = \cos(\pi(s + u_1 +v +\frac{w}{2}))
\Gamma(s + u_1 + v + \frac{w}{2} + it) \Gamma(s + u_1 + v + \frac{w}{2} - it)
\\
\Gamma(\tfrac12 - s -u_1 -v)^2 
\frac{\Gamma(v) \Gamma(\half + s + u_2 -v - \frac{w}{2})}{\Gamma(\half + s + u_2  - \frac{w}{2})}
\frac{G(s) G(w)}{s w} g(s) \widetilde{W}(u_1, u_2) \zeta(1 + w)
c^{*},
\end{multline*}
and where $c^*$ is meant to account for bounded factors like powers of $2$, $\pi$, etc. that do not have any effect on the convergence of the integrals.  

Similarly, we compute
\begin{multline*}
K_{r-}(t) =  \left(\frac{1}{2 \pi i} \right)^5 \int_{(c_s)} \int_{(c_w)} \int_{(c_{u_2})} \int_{(c_{u_1})} \int_{(c_v)}    \frac{q^{2s} M^{u_1} N^{u_2} H_{-}(s, u_1, u_2, v, w)}{(4\pi)^{2s +2u_1 + 2v + w} f^{2s +u_1 + u_2} m^{\frac{w}{2}}}  
\\
2 \int_0^{\infty}K_{2it}(x) 
x^{2s +2u_1 +  2v +w} \frac{dx}{x}
dv du_1 du_2 dw ds.
\end{multline*}
This time we use the formula (see 6.561.16 of \cite{GR})
\begin{equation*}
2\int_0^{\infty} K_{2it}(x) x^{s} \frac{dx}{x} = 2^{s-1} \Gamma\left(\frac{s + 2it}{2}\right)\Gamma\left(\frac{s - 2it}{2}\right),
\end{equation*}
valid for $\text{Re}(s) > |\text{Re}(2it)|=0$.  Hence
\begin{equation}
\label{eq:Ktransform}
K_{r-}(t) =  \left(\frac{1}{2 \pi i} \right)^5 
\int \dots \int
\frac{q^{2s} M^{u_1} N^{u_2}}{f^{2s +u_1 + u_2} m^{\frac{w}{2}}} \widehat{H}_{-}(s, u_1, u_2, v, w;t)
dv du_1 du_2 dw ds,
\end{equation}
where
\begin{multline*}
\widehat{H}_{-}(s,u,v,w; t) = 
\Gamma(s + u_1 + v + \frac{w}{2} + it) \Gamma(s + u_1 + v + \frac{w}{2} - it)
\Gamma(\tfrac12 - s -u_1 -v)^2 \sin(\pi(s + u_1 +v))
\\
\frac{\Gamma(v) \Gamma(\half + s + u_2 -v - \frac{w}{2})}{\Gamma(\half + s + u_2  - \frac{w}{2})}
\frac{G(s) G(w)}{s w} g(s) \widetilde{W}(u_1, u_2) \zeta(1 + w)
c^{*}.
\end{multline*}

\subsection{Maass forms}
This section is devoted to proving
\begin{myprop}  We have
\begin{equation}
E_{m\pm} \ll q^{-\half + \theta + \varepsilon} \left(\frac{N}{M}\right)^{\half},
\end{equation}
\end{myprop}
We treat the opposite sign case $E_{m-}$ only since the case of $E_{m+}$ is similar, and easier.  
\begin{proof}
We have
\begin{equation*}
E_{m-} = \frac{1}{\phi^{*}(q)}  \sum_{d | q} \phi(d) \mu(\frac{q}{d})
\sum_{f \equiv 0 \shortmod{d}} \sum_m \frac{d(m)}{m^{\half} f^{\half}}
\sum_j \rho_{j}(m) \overline{\rho_{j}(-f)} K_{r-}(\kappa_j),
\end{equation*}
which upon using \eqref{eq:Ktransform} is
\begin{multline}
\label{eq:scanner}
E_{m-} =  \frac{1}{\phi^{*}(q)}  \sum_{d | q} \phi(d) \mu(\frac{q}{d})
\sum_{f \equiv 0 \shortmod{d}} \sum_m \frac{d(m)}{m^{\half} f^{\half}}
\sum_j  \frac{|\rho_{j}(1)|^2 \lambda_j(m) \lambda_j(f)}{\cosh{\pi \kappa_j}} 
\\
 \left(\frac{1}{2 \pi i} \right)^5 
\int \dots \int
\frac{q^{2s} M^{u_1} N^{u_2}}{f^{2s +u_1 + u_2} m^{\frac{w}{2}}}    \widehat{H}_{-}(s, u_1, u_2, v, w;\kappa_j) \cosh(\pi \kappa_j) 
dv du_1 du_2 dw ds.
\end{multline}
Let $E_{K}$ be the same expression as \eqref{eq:scanner} but with the spectral parameter $\kappa_j$ restricted to the dyadic segment $K \leq \kappa_j < 2K$.  Taking $c_s = \frac38$ and $c_w = \frac32$ means the sums over $m$ and $f$ converge absolutely.  Now we use the following variation on \eqref{eq:Ramanujan}
\begin{equation*}
\sum_n \frac{\sigma_{\lambda}(n) \lambda_j(n)}{n^s} = \frac{L_j(s) L_j(s - \lambda)}{\zeta(2s-\lambda)}
\end{equation*}
and
\begin{multline*}
\sum_{f \equiv 0 \shortmod{q}} \frac{\lambda_j(f)}{f^s} = q^{-s} \Big(\sum_{(f,q) = 1}  \frac{\lambda_j(f)}{f^s} \Big)
\Big(\sum_{n \geq 0}  \frac{\lambda_j(q^{n+1})}{q^{ns}} \Big)
\\
=q^{-s} L_{j, q}(s) \sum_{n \geq 0}  \frac{\lambda_j(q) \lambda_j(q^{n}) - \lambda_j(q^{n-1})}{q^{ns}}
=
q^{-s} L_{j}(s) \left( \lambda_j(q)  - q^{-s} \right),
\end{multline*}
to obtain
\begin{multline*}
E_{K} =  \frac{1}{\phi^{*}(q)}  \sum_{d | q} \frac{\phi(d)}{d^{\half}} \mu(\frac{q}{d})
\sum_{K \leq \kappa_j < 2K}  \frac{|\rho_{j}(1)|^2}{\cosh{\pi \kappa_j}}
\left(\frac{1}{2 \pi i} \right)^5 \int_{(c_s)} \int_{(c_w)} 
\\
\int_{(c_{u_2})} \int_{(c_{u_1})} \int_{(c_v)} 
\frac{q^{2s}M^{u_1} N^{u_2}}{d^{2s + u_1 + u_2}}  
L_j(\tfrac12 + \tfrac{w}{2})^2 L_j(\tfrac12+ 2s + u_1 + u_2) 
\\
\Big(\lambda_j(d) -\frac{\delta_{d,q}}{d^{\half + 2s + u_1 + u_2}} \Big) 
\frac{ \widehat{H}_{-}(s, u_1, u_2, v, w;\kappa_j) \cosh(\pi \kappa_j)}{\zeta(1 + w)} 
dv du_1 du_2 dw ds.
\end{multline*}
To estimate this term we initially move the lines of integration so that 
\begin{equation}
\label{eq:linesatc}
c_s = c_v= c_{w} = \varepsilon, \quad c_{u_1} =  c, \quad c_{u_2} = -c+\varepsilon,
\end{equation}
where $3\epsilon < c < \half-3\epsilon$,
passing no poles in this process by close examination of the form of $\widehat{H}_{-}$.  Due to the rapid decay of $\widehat{H}_{-}$, we may truncate the integrals so that $\text{Im}(s)$, $\text{Im}(u_1)$, $\text{Im}(u_2)$, $\text{Im}(v)$, $\text{Im}(w) \ll (qK)^{\varepsilon}$ with a negligible error (say, size $\ll (qK)^{-1000}$).  The issue at hand is the dependence on $K$.  
Using Stirling's approximation, we see that
\begin{equation}
\label{eq:stirlingsapproximation}
\cosh(\pi \kappa_j) \Gamma(s + u_1 + v + \tfrac{w}{2} + i\kappa_j) \Gamma(s + u_1 + v + \tfrac{w}{2} - i\kappa_j) \ll 
q^{\varepsilon} K^{-1 + 2 c}.
\end{equation}
Then we have
\begin{equation}
\label{eq:thirdmomentformula}
E_K \ll q^{-\half} (qK)^{\varepsilon} K^{-1 + 2 c} \left(\frac{N}{M}\right)^{-c} \max_{s_1, s_2} \Big| \sum_{K \leq \kappa_j < 2K} \frac{|\rho_{j}(1)|^2}{\cosh{\pi \kappa_j}} \lambda_j(q) L_j(\tfrac12 + s_1) L_j(\tfrac12 + s_2)^2 \Big|,
\end{equation}
where the maximum over $s_1$ and $s_2$ is over the rectangles $0 \leq \text{Re}(s_i) \leq \varepsilon$ and $|\text{Im}(s_i)| \ll (qK)^{\varepsilon}$, $i=1,2$.  

The spectral sum may be estimated using $\lambda_j(q) \ll q^{\theta + \varepsilon}$ and the fourth moment bound
\begin{equation}
\label{eq:maassfourthestimate}
\sum_{K \leq \kappa_j \leq 2K} \frac{|\rho_{j}(1)|^2}{\cosh{\pi \kappa_j}}  |L_j(\tfrac12 + s)|^4 \ll K^{2 + \varepsilon},
\end{equation}
for $\text{Im}(s) \ll K^{\varepsilon}$, $\text{Re}(s) \geq 0$, which follows from the large sieve inequality for Maass forms \cite{IwaniecMean}.  For a proof, see Theorem 3.4 in Motohashi's book \cite{M} (actually Motohashi had $s=0$ but this is a minor change).

First suppose that $K \gg M^{-\half} N^{\half} q^{\varepsilon}$.  Observation of the location of poles $\widehat{H}_{-}$ shows that we can take $c = A$ large without crossing any poles.
For such $K$ we thus obtain
\begin{equation}
E_{K} \ll (qK)^{-1000}.
\end{equation}
For the complementary range, taking $c = \half - 3\varepsilon$ completes the proof with the bound
\begin{equation}
E_{K} \ll q^{-\half + \theta}  K^2 \left(\frac{N}{M}\right)^{-\half} (qK)^{\varepsilon} \ll q^{-\half + \theta + \varepsilon} \left(\frac{N}{M}\right)^{\half}. \qedhere
\end{equation}
\end{proof}
The estimate of $E_{m+}$ is even easier than that for $E_{m-}$ because $\cosh(\pi t) M_{r+}(t)$ has exponential decay as $t \rightarrow \infty$.

\subsection{The continuous spectrum}
In this section we prove
\begin{myprop}  We have
\begin{equation}
E_{c\pm} \ll q^{-\half + \varepsilon} \left(\frac{N}{M}\right)^{\frac14}.
\end{equation}
\end{myprop}
The exponent $\frac14$ depends on an estimate for the $6$th moment of the Riemann zeta function.
As in the previous section, we only show the full details for $E_{c-}$ since $E_{c+}$ is similar and even easier to handle.
\begin{proof}
Our starting point is
\begin{equation*}
E_{c-} = \frac{1}{\phi^{*}(q)}  \sum_{d | q} \phi(d) \mu(\frac{q}{d})
\sum_{f \equiv 0 \shortmod{d}} \sum_m \frac{d(m)}{m^{\half} f^{\half}}
\frac{1}{\pi} \intR \frac{\sigma_{2it}(m) \sigma_{2it}(f)}{(mf)^{it} |\zeta(1 + 2it)|^2} \cosh(\pi t)  K_{r-}(t) dt,
\end{equation*}
which after using \eqref{eq:Ktransform} becomes
\begin{multline*}
E_{c-} = \frac{1}{\phi^{*}(q)}  \sum_{d | q} \phi(d) \mu(\frac{q}{d})
\sum_{f \equiv 0 \shortmod{d}} \sum_m \frac{d(m)}{m^{\half} f^{\half}}
\frac{1}{\pi} \intR \frac{\sigma_{2it}(m) \sigma_{2it}(f)}{(mf)^{it} |\zeta(1 + 2it)|^2}
\\
\left(\frac{1}{2 \pi i} \right)^5 
\int \dots \int
\frac{q^{2s} M^{u_1} N^{u_2} }{f^{2s +u_1 + u_2} m^{\frac{w}{2}}}  \widehat{H}_{-}(s, u_1, u_2, v, w;t) \cosh(\pi t)
dv du_1 du_2 dw ds dt.
\end{multline*}
Now move $c_s$ to $\half$ and $c_w$ to $1 + \varepsilon$ and execute the summations over $m$ and $f$ in terms of a product of zeta functions to get
\begin{multline*}
E_{c-} = \frac{1}{\phi^{*}(q)}  \sum_{d | q} \frac{\phi(d)}{d^{\half}} \mu(\frac{q}{d})
\frac{1}{\pi} \intR 
\left(\frac{1}{2 \pi i} \right)^5 \int_{(c_s)} \int_{(c_w)} \int_{(c_{u_2})} \int_{(c_{u_1})} \int_{(c_v)} 
\\
\frac{q^{2s} M^{u_1} N^{u_2}}{d^{2s + u_1 + u_2}} \widehat{H}_{-}(s, u_1, u_2, v, w;t) \cosh(\pi t) \frac{Z(s,u_1, u_2, w,t)}{|\zeta(1 + 2it)|^{2}}
dv du_1 du_2 dw ds dt,
\end{multline*}
where
\begin{multline*}
Z(s, u_1, u_2, w, t) = \zeta({\textstyle \half + \frac{w}{2} + it})^2
\zeta({\textstyle \half + \frac{w}{2} - it})^2
\\
\zeta({\textstyle \half + 2s + u_1 + u_2 + it})\zeta({\textstyle \half + 2s + u_1 + u_2 - it}) A_d(s,u_1,u_2,w),
\end{multline*}
where $A_d(s,u_1,u_2,w)$ is bounded by $d^{\varepsilon}$ and holomorphic for $\text{Re}(s), \text{Re}(u_1),\text{Re}(u_2), \text{Re}(w) > - \varepsilon$ (use Lemma \ref{lem:divisorfunctionsum} to get an exact expression).  As in the previous section, take the contours according to \eqref{eq:linesatc}
where $3\varepsilon < c < \half-3\varepsilon$ is at our disposal.  Although we crossed various poles of the zeta function, since they are all at height $t$ (roughly speaking), the decay of the test functions shows that these terms have rapid decay in $t$-aspect, and it is not difficult to bound the contribution of these terms by $O(q^{-1000})$.

Again, we may truncate all the integrals except the one over $t$ at height $(qt)^{\varepsilon}$ with negligible error.
The crucial issue is convergence in $t$-aspect.
Suppose we have a bound on the $6$th moment of the Riemann zeta function of the form
\begin{equation}
\frac{1}{T} \int_T^{2T} |\zeta(\tfrac12 + it)|^6 dt \ll T^{\theta' + \varepsilon}.
\end{equation}
The currently best-known result has $\theta' = 1/4$ using H\"{o}lder's inequality with
\begin{align*}
\int_T^{2T} |\zeta(\tfrac12 + it)|^4 dt \ll T^{1+ \varepsilon}, \quad
\int_T^{2T} |\zeta(\tfrac12 + it)|^{12} dt \ll T^{2+ \varepsilon},
\end{align*}
This bound on the $12$th moment was proved by Heath-Brown \cite{H-B3}.

Using \eqref{eq:stirlingsapproximation} shows that $c=-\theta' - \varepsilon$ is sufficient for convergence in $t$-aspect, and gives the desired bound.  If $\theta' = 0$ (a standard conjecture) it would give the bound $E_{c\pm} \ll q^{-\half + \varepsilon}$.
\end{proof}

\subsection{The dual terms}
In this section we sketch how to treat $E_{\overline{M,N}}$.  The basic outline is the same as for $E_{M,N}$ but some convergence issues are slightly more delicate and the arguments must be modified in some places.

Using the same computations as in Section \ref{section:reduction}, we get that the analog of \eqref{eq:theonion} is
\begin{multline*}
\overline{E}_{\pm} = \frac{1}{\phi^{*}(q)}  \sum_{d | q} \phi(d) \mu(\frac{q}{d}) \sum_{l} \frac{1}{l}  
\sum_{\substack{f \equiv 0 \shortmod{d}}}{} \sum_m \frac{S(m,\mp f;l) }{f^{\half} m^{\half}}
\\
\left(\frac{1}{2 \pi i} \right)^5 \int_{(c_s)} \int_{(c_w)} \int_{(c_{u_2})} \int_{(c_{u_1})} \int_{(c_v)}   \frac{q^{2s} M^{u_1} N^{u_2} m^{s + u_1 + v}}{f^{s+ u_2 - v - \frac{w}{2} } l^{2s + 2u_1 + 2v +w}  }
\\
\overline{H}_{\pm}(s, u_1, u_2, v, w)
dv du_1 du_2 dw ds,
\end{multline*}
where $\overline{H}_{\pm}$ is given by an expression identical to \eqref{eq:Hpmdefinition} except $H$ is replaced by $\overline{H}$.  More explicitly, the ratio of gamma factors
\begin{equation}
\label{eq:gammaplus}
\frac{\Gamma(v) \Gamma(\half + s + u_2 -v - \frac{w}{2})}{\Gamma(\half + s + u_2  - \frac{w}{2})}
\end{equation}
appearing in the definition of $H_{\pm}$ is replaced by
\begin{equation}
\label{eq:gammaminus}
\frac{\Gamma(v) \Gamma(\half - s - u_2 + \frac{w}{2})}{\Gamma(\half - s - u_2 +v + \frac{w}{2})}.
\end{equation}
Furthermore, the term $S(m,f;l)$ is matched with $S_{-}$ and $S(m,-f;l)$ is matched with $S_{+}$.  The primary difference between \eqref{eq:gammaplus} and \eqref{eq:gammaminus} is that the former expression has exponential decay in $v$-aspect, and the latter has at best polynomial decay.

Because the opposite sign case is matched with $S_{+}$, the convergence in $v$-aspect is assured by the presence of the factor $\Gamma(\half -s - u_1 - v)^2$, and the same methods as before work to give the bound \eqref{eq:Ebound}.

The same sign case reduces to
\begin{multline*}
\frac{1}{\phi^{*}(q)}  \sum_{d | q} \frac{\phi(d)}{d^{\half}} \mu(\frac{q}{d})
\sum_{K \leq \kappa_j < 2K}  \frac{|\rho_{j}(1)|^2}{\cosh{\pi \kappa_j}}
\\ 
\left(\frac{1}{2 \pi i} \right)^5 \int_{(c_s)} \int_{(c_w)} \int_{(c_{u_2})} \int_{(c_{u_1})} \int_{(c_v)} 
\frac{q^{2s}M^{u_1} N^{u_2}}{d^{2s + u_1 + u_2}}  
L_j(\tfrac12 + \tfrac{w}{2})^2 L_j(\tfrac12 + 2s + u_1 + u_2) 
\\
(\lambda_j(d) -\frac{\delta_{d,q}}{d^{\half + 2s + u_1 + u_2}}) 
\frac{ \widehat{\overline{H}}_{-}(s, u_1, u_2, v, w;\kappa_j) \cosh(\pi \kappa_j)}{\zeta(1 + w)} 
dv du_1 du_2 dw ds,
\end{multline*}
where
\begin{multline*}
\widehat{\overline{H}}_{-}(s,u,v,w; t) = \cos(\pi(s + u_1 +v +\frac{w}{2}))
\Gamma(s + u_1 + v + \frac{w}{2} + it) \Gamma(s + u_1 + v + \frac{w}{2} - it)
\\
\sin(\pi(s + u_1 + v))
\Gamma(\tfrac12 - s -u_1 -v)^2 
\frac{\Gamma(v) \Gamma(\tfrac12 - s - u_2 + \frac{w}{2})}{\Gamma(\tfrac12 - s - u_2  +v + \frac{w}{2})}
\frac{G(s) G(w)}{s w} g(s) \widetilde{W}(u_1, u_2) \zeta(1 + w)
c^{*}.
\end{multline*}
Again, we truncate the integrals over $s, u_1, u_2, w$ at imaginary part $\ll (qK)^{\varepsilon}$.  The convergence in $v$-aspect can be detected via Stirling's approximation, which gives
\begin{multline*}
\cos(\pi(s + u_1 +v +\frac{w}{2}))
\Gamma(s + u_1 + v + \frac{w}{2} + it) \Gamma(s + u_1 + v + \frac{w}{2} - it)
\\
\sin(\pi(s + u_1 + v))
\Gamma(\tfrac12 - s -u_1 -v)^2 
\frac{\Gamma(v) \Gamma(\half - s - u_2 + \frac{w}{2})}{\Gamma(\half - s - u_2  +v + \frac{w}{2})}
\\
\ll
(qK)^{\varepsilon} e^{\pi |v|} e^{-\frac{\pi}{2} |v - t|} e^{-\frac{\pi}{2} |v + t|}
(1 + |v|)^{-\frac{3}{2} + c_{u_2}}.
\end{multline*}
A careful but elementary estimation gives the bound
\begin{equation*}
\int_0^{\infty} e^{\pi y} e^{-\frac{\pi}{2} |y - t|} e^{-\frac{\pi}{2} |y + t|}
(1 + y)^{-\frac{3}{2} + c_{u_2}} dy \ll (1 + t)^{-\frac{3}{2} + c_{u_2}}.
\end{equation*}
Thus we may take $c_{u_2} = -\half - 3\varepsilon$ and $c_{u_1} = \half - 3\varepsilon$ to get
\begin{equation*}
E_{\overline{M,N}} \ll q^{-\half + \theta + \varepsilon} M^{\half} N^{-\half}.
\end{equation*}
Notice this is better than \eqref{eq:Ebound}.  Actually, this phenomenon holds true for $E_{m+}$ also: for this term, the convergence in $\kappa_j$ aspect holds even when all the lines of integration are at $\varepsilon$, say, and the loss of $(N/M)^{\half}$ does not appear in the bound.

The treatment of the continuous spectrum follows similar lines.


\end{document}